\documentclass[twoside, 11pt]{article}
\usepackage[utf8]{inputenc}
\usepackage[T1]{fontenc}
\usepackage[charter]{mathdesign}
\usepackage{amsthm, amsfonts, enumerate, calligra,mathrsfs}
\usepackage[a4paper,margin=3.4cm, headsep=15pt, headheight=2cm]{geometry}
\newcommand\smvee{\raise0.3ex\hbox{$\scriptscriptstyle\vee$}}
\usepackage{fancyhdr}
\usepackage{tikz-cd, tikz}
\usetikzlibrary{matrix, calc}
\usepackage{mathtools}

\newtheorem{theorem}{Theorem}[section]
\newtheorem{nmtheorem}{Theorem}

\newtheorem{prop}[theorem]{Proposition}
\newtheorem{remark}[theorem]{Remark}
\newtheorem{lemma}[theorem]{Lemma}

\newtheorem{cor}[nmtheorem]{Corollary}
\newtheorem{defn}[theorem]{Definition}

\numberwithin{equation}{section}
\newcounter{para}

\usepackage[colorlinks,backref]{hyperref}
\usepackage[nameinlink]{cleveref}
\hypersetup{colorlinks,breaklinks,
	citecolor=[rgb]{0.1,0.7,0.5},
	linkcolor=[rgb]{0.3,0.4,0.7}}
\newcommand{\etale}{étale }

\def\C{\mathbb{C}}
\def\Q{\mathbf{Q}}
\def\Qb{\mathbb{Q}}
\def\A{\mathbb{A}}
\def\fS{\mathbb{S}}
\def\grd{\mathfrak{g}^r_d}
\def\P{\mathbb{P}}

\def\N{\mathbb{N}}

\def\V{\mathcal{V}}
\def\F{\mathcal{F}}

\def\Deltainj{\Delta_{\mathrm{inj}}}
\def\et{\acute{e}t}

\def\1{\mathbb{1}}

\def\txf{t_x^\mathit{eff}}
\def\morn{\mathrm{Mor}_n(C,\P^r)}

\def\sgn{\mathit{sgn}}
\def\pic{\mathit{Pic}}
\def\deg{\mathit{deg\,}}

\def\uconf{\mathit{UConf}}
\def\Aut{\mathit{Aut}}
\newcommand\Sym{\mathrm{Sym\,}}

\date{}
\DeclareMathOperator{\Shom}{\mathscr{H}\text{\kern -3pt {\calligra\large om}}\,}
\DeclareMathOperator{\Sext}{\mathscr{E}\text{\kern -3pt {\calligra\large xt}}\,}
\newcount \questionno
\def\question{\medbreak
	\global \advance \questionno 1
	\textbf{Problem \the\questionno}.\enspace \ignorespaces}
\newcommand\shorttitle{Filtration of cohomology via symmetric semisimplicial spaces}
\newcommand\authors{Oishee Banerjee}
\usepackage{titlesec}

\titleformat{\section}
{\normalfont\fontsize{11}{13}\bfseries}{\thesection}{1em}{}
\title{Filtration of cohomology via symmetric semisimplicial spaces}
\author{Oishee Banerjee}

\begin{document}
	
	\maketitle
	\begin{abstract}
		In the simplicial theory of hypercoverings we replace the indexing category $\Delta$ by the \emph{symmetric simplicial category} $\Delta S$ and study (a class of) $\Deltainj S$-hypercoverings, which we call \emph{spaces admitting symmetric (semi)simplicial filtration}- this special class happens to have a structure of a module over a graded commutative monoid of the form $\Sym M$ for some space $M$. For $\Delta S$-hypercoverings we construct a spectral sequence, somewhat like the Čech-to-derived category spectral sequence. The advantage of working with $\Delta S$ over $\Delta$ is that various combinatorial complexities that come with working on $\Delta$ are bypassed, giving simpler, unified proof of results like the computation of (in some cases, stable) singular cohomology (with $\Qb$ coefficients) and \etale cohomology (with $\Qb_{\ell}$ coefficients) of the moduli space of degree $n$ maps $C\to \P^r$ with $C$ a smooth projective curve of genus $g$, of unordered configuration spaces, of the moduli space of smooth sections of a fixed $\grd$ that is $m$-very ample for some $m$ etc.  In the special case when a $\Deltainj S$-object  $X$ \emph{admits a symmetric semisimplicial filtration by $M$}, we relate these moduli spaces to a certain derived tensor.
	\end{abstract}
	
	\tableofcontents

	\section{Introduction}\label{intro}
	
	The theory of simplicial spaces forms the core of Verdier's theory of hypercoverings and the subsequent vast generalisations in Deligne's theory of cohomological descent.  We build a theory of hypercoverings where the traditional indexing category $\Delta$ (commonly known as the simplicial category) is replaced by $\Delta S$ (which we call the \emph{symmetric simplicial category} (see Definition \ref{def2.1}), and which contains $\Delta$ as a subcategory). Intuitively speaking, whereas $\Delta S$ enjoys all the properties that makes $\Delta$ a fundamental part of homotopy theory, its objects also have nontrivial automorphism groups, isomorphic to the symmetric groups.  \footnote{This might give the reader the impression that $\Delta S$ is the equivalent to the category of finite sets, but it's not, as one can gather immediately from the axioms of Definition \ref{def2.1}}. And this is the main advantage of working with $\Delta S$. The category $\Delta S$ was first introduced as a part of the concept of \emph{crossed simplicial groups} independently by Fiedorowicz-Loday (\cite{FiedorowiczLoday1991}) and Krasauskas (\cite{Krasauskas.1987}). 
	In particular, if one's goal is, for example, to compute (stable) (co)homology of moduli spaces (which often come in families indexed by a parameter, say $n$) that are naturally quotients of spaces equipped with permutation actions by $\{S_n\}$, then $\Delta S$, by encoding the permutation action as automorphisms in the category itself, gives us a precise tool to entirely bypass all the combinatorial complexities that form a part of $\Delta$. 
	
	Amongst $\Delta S$ spaces, we define, and give special attention to, \emph{spaces admitting symmetric semisimplicial filtration} (see Definition \ref{def2.10}) because of their frequent manifestations in topology and geometry. Roughly speaking, given a family of spaces $\{X_n\}_{n\in \N}$, we say $X:= \sqcup X_n$ (or equivalently $\{X_n\}$)admits symmetric semisimplicial filtration by a space $M$ with filter gap $e$, if $X$ forms a module over the graded commutative topological monoid $\Sym M$, where $M$ has grading $e$, and satisfies two minor additional conditions (see Definition \ref{def2.10}). 
	We call $$U_n:= X_n - f_{0}(M\times X_{n-e})$$ the space of $M$-indecomposables of $X_n$; we use the same term for $$U:=\sqcup U_n$$ as well. We also write $$U_{} = X_{}-(M\times X_{})$$ which has the unambiguous meaning of $ X_{} - f_{0}(M\times X_{})$ under the module structure.
	
	\paragraph{Results.} 
	To state our first two theorems we need some notations and conventions.
	For a graded vector space $V$, let $V^{(r)}$ denote its $r^{th}$ graded component,  and let $V^{\textrm{odd}}: = \oplus_{j\in \mathbb{Z}} V^{(2j+1)}$ and $V^{\textrm{even}}:= \oplus_{j\in \mathbb{Z}} V^{(2j)}$ denote the odd and even graded subspaces of $V$, respectively. Throughout this paper, by a \emph{space} we mean a locally-compact Hausdorff topological space or a quasi-projective algebraic variety over some field. By a \emph{morphism} we mean a continuous map of topological spaces or a morphism of algebraic varieties. For a $\mathbb{Z}$-scheme $X$ we continue to denote its base change to any algebraically closed field $K$ by $X$; in turn we mean $H^q(X;\Q)$ (respectively,  $H^q_c(X;\Q)$) to stand for both the singular cohomology $H^q(X(\mathbb{C});\mathbb{Q})$ (respectively, $H^q_c(X(\mathbb{C});\mathbb{Q})$, singular cohomology with compact support) as well as the \etale cohomology $H_{\et}^q(X(K); \mathbb{Q}_{\ell})$ (respectively, $H_{\et,c}^q(X(K); \mathbb{Q}_{\ell})$, \etale cohomology with compact support),  $\ell$ coprime to $\mathrm{char}\,\, K$. 
	
	\begin{nmtheorem}[\textbf{Cohomology of indecomposables vs. indecomposables in cohomology}]\label{thm1}
		Let $M$ and $\{X_n\}_{n\in\mathbb{N}}$ be locally compact connected Hausdorff topological spaces and let $X=\sqcup X_n$. Suppose that $X$ admits a semisimplicial filtration by $M$, with face maps given by $$f_i:M^p\times X_{n-ep} \to M^{p-1}\times X_{n-e(p-1)}.$$ Let $e>0$ be the filter gap and $U=\sqcup U_n$ the space of $M$-indecomposables. Then there exists a spectral sequence \begin{gather}\label{eq1.1}
			E_1^{p,q} =\bigoplus_{l+m=q}\,\,\bigoplus_{i+j=p}  \big(\Sym^i H_c^{\textrm{odd}}(M;\mathbb{Q}) \otimes \Lambda^j H_c^{\textrm{even}}(M;\mathbb{Q})\big)^{(l)} \otimes H^m_c(X_{n-ep};\mathbb{Q}) \implies H^{p+q}_c(U_n;\mathbb{Q})
		\end{gather} where the differentials are given by alternating sum of the pullbacks on cohomology induced by the face maps: \begin{gather*}
			d_1^{p,q}: E_1^{p,q} \to E_1^{p+1,q} \\ d_1^{p,q}:= \sum_{i=0}^{p-1}(-1)^i f_i^*.
		\end{gather*}
		If $\{X_n\}$ and $M$ are quasi-projective algebraic varieties over a field $K$, then there is a spectral sequence of $Gal(\overline{K}/K)$-representations \begin{gather*}
			E_1^{p,q} =\bigoplus_{l+m=q}\,\,\bigoplus_{i+j=p}  \big(\Sym^i H_{\et,c}^{\textrm{odd}}(M;\mathbb{Q}_{\ell}) \otimes \Lambda^j H_{\et,c}^{\textrm{even}}(M;\mathbb{Q}_{\ell})\big)^{(l)} \otimes H^m_{\et,c}(X_{n-ep};\mathbb{Q}_{\ell}) \\ \implies H^{p+q}_{\et,c}(U_n;\mathbb{Q_{\ell}}) 
		\end{gather*}where $\ell$ is coprime to $\mathrm{char}\, K$, and the differentials are exactly the same as above.
		$\hfill\square$
	\end{nmtheorem}
	In the special case when all spaces are smooth projective varieties or compact oriented manifolds without boundaries, one obtains  close cousin (essentially the Verdier dual) of Theorem \ref{thm1} as follows.
	
	\begin{nmtheorem}\label{thm2}
		Let $M$ and $\{X_n\}_{n\in\mathbb{N}}$ be compact oriented manifolds without boundaries. Suppose that $\{X_n\}_{n\in\mathbb{N}}$ admits a semisimplicial filtration by $M$, with face maps given by $$f_i:M^p\times X_{n-ep} \to M^{p-1}\times X_{n-e(p-1)}.$$ Let $e>0$ be the filter gap and $\{U_n\}$ the space of $M$-indecomposables. Furthermore, let $$c(n,p):= \dim_{\mathbb{R}}(X_n) - \dim_{\mathbb{R}} (M^p\times X_{n-ep}).$$ Then there exists a second quadrant spectral sequence which converges to $H^*(U_n;\Qb)$ as an algebra. The $E_1$ page of that spectral sequence reads as:
		
		\begin{gather}\label{eq1.2}
			E_1^{-p,q} =\bigoplus_{l+m=q-c(n,p)}\,\,\bigoplus_{i+j=p}  \big(\Sym^i H^{\textrm{odd}}(M;\mathbb{Q}) \otimes \Lambda^j H^{\textrm{even}}(M;\mathbb{Q})\big)^{(l)} \otimes H^m(X_{n-ep};\mathbb{Q}) \nonumber \\ \implies H^{q+p}(U_n;\Qb)
		\end{gather} with the differentials given by the alternating sum of the Gysin pushforwards induced by the face maps i.e. 
		\begin{gather*}
			d_1^{-p,q}: E_1^{-p,q} \to E_1^{-(p-1),q} \\ d_1^{p,q}:= \sum_{i=0}^{p-1}(-1)^i {f_i}_*.
		\end{gather*}
		
		If $\{X_n\}$ and $M$ are smooth projective algebraic varieties over a field $K$, and if we define $$c(n,p):= \dim_{\overline{K}}(X_n) - \dim_{\overline{K}} (M^{p}\times X_{n-ep}),$$ then we have a second quadrant spectral sequence of $Gal(\overline{K}/K)$-representations whose $E_1$ page reads as \begin{gather*}
			E_1^{-p,q} =\bigoplus_{l+m=q-2c(n,p)}\,\,\bigoplus_{i+j=p}  \big(\Sym^i H_{\et}^{\textrm{odd}}(M;\mathbb{Q}_{\ell}) \otimes \Lambda^j H_{\et}^{\textrm{even}}(M;\mathbb{Q}_{\ell})\big)^{(l)} \otimes H^m_{\et}(X_{n-ep};\mathbb{Q}_{\ell}) (-c(n,p)) \\ \implies H^{q+p}_{\et}(U_n;\mathbb{Q_{\ell}}) 
		\end{gather*}where $\ell$ is coprime to $\mathrm{char}\, K$, the differentials exactly the same as above, and the spectral sequence converges to $H^*(U_n;\Qb_{\ell})$ as an algebra.
		$\hfill\square$
	\end{nmtheorem}
	\begin{remark}
		
		To the reader familiar with the concept of derived indecomposables (see, for example, \cite[Defnition 8.5]{GKRW2021}) note that the dual (in the sense of Verdier duality) of the spectral sequence above is $$\mathrm{Tor}^{H^{\mathrm{BM}}_*(\Sym M)}(\Q,\Q) \otimes  H^{\mathrm{BM}}_*(X) \cong \Sym(H^{\mathrm{BM}}_*(\Sigma M))\otimes H^{\mathrm{BM}}_*(X)$$ where  $\Sigma M$ denotes the suspension of $M$; which is the associated graded of $$\mathrm{Tor}^{H^{\mathrm{BM}}_*(\Sym M)}(H^{\mathrm{BM}}_*(X), \Q),$$ the derived indecomposables of $H^{\mathrm{BM}}_*(X)$ as a $H^{\mathrm{BM}}_*(\Sym M)$-module. It is not true in general that $\mathrm{Tor}^{H^{\mathrm{BM}}_*(\Sym M)}(H^{\mathrm{BM}}_*(X), \Q) \cong H^{\mathrm{BM}}_*(U, \Q)$- a similar isomorphism does, however, hold if we replace homology by the corresponding chain complexes: $$\mathrm{Tor}^{\Sym(\mathrm{R}\Gamma_c(M,\omega_M))}(\mathrm{R}\Gamma_c(X,\omega_X), \Q) \cong \mathrm{R}\Gamma_c(U,\omega_U)$$ where, for any space $B$, $\omega_B$ denotes its dualizing sheaf. This isomorphism follows from Lemma \ref{lem2.11}, and is proved in a more general setting in \cite{Banerjee23}.
	\end{remark}
	
	
	Before we state the other results, let us briefly look at the ubiquity of families that admit a symmetric semisimplicial filtration. 
	\paragraph{Some context and some examples.} There are many examples of families of spaces admitting a symmetric semisimplicial filtration (and thus satisfying the hypothesis of Theorem \ref{thm1}), including, but not limited to the following.
	\begin{enumerate}
		
		\item 	\emph{The $n^{th}$-symmetric powers of a space $X$.} Let $X_n = \Sym^nM$. Define \begin{gather}
			f_{i}: M^{p+1}\times \Sym^{n-2(p+1)}M \to M^{p} \times \Sym^{n-2p}M \nonumber \\\  (a_0,\ldots, a_p) , \{b_1,\ldots, b_{n-2p}\} \mapsto  (a_0,\ldots, \hat{a_i}, \ldots, a_p), \{ a_i, a_i , b_1,\ldots, b_{n-2p}\}\label{eq1.3}
		\end{gather} where $(a_0,\ldots, a_p)$ denotes an ordered $p+1$-tuple of elements in $M$, and $\{b_1,\ldots, b_p\}$  denotes an unordered $p$-tuple and $\hat{a_i}$ stands for $a_i$ (the $(i+1)^{th}$ entry) removed. It is easy to check that with these morphisms as face maps, the semisimplicial space $\{M^p\times \Sym^{n-2p}M\} $ naturally conforms to Definition \ref{def2.10} i.e. $\Sym M$ admits symmetric simplicial filtration by $M$ with filter gap $2$. The space of $M$-indecomposables is $U_n = \uconf_n(M)$, the unordered configuration space of $n$ distinct points in $M$. For the explicit computation of the spectral sequence that converges to $H^*(\uconf_n(M);\Q)$, see Corollary \ref{cor6}.
		
		\item \emph{The moduli space of $(r+1)$-tuples of monic polynomials of degree $n$.} Let $\mathit{Poly}^{n,r+1}$ be the space of $(r+1)$-tuples of monic degree $n$ homogeneous polynomials in one variable over an algebraically closed field $K$, and let $\mathit{Poly}_v^{n,r+1}$ be the locus of those $r$-tuples having no common roots of multiplicity $\geq v$. Then $\mathit{Poly}^{n,r+1}$ admits a symmetric semisimplicial filtration by $\A^1$. Indeed, we have a semisimplicial space given by $\{(\A^1)^p\times \mathit{Poly}^{n-pv,r+1}\}_{0\leq p\leq n}$ with face maps defined by \begin{gather*}
			f_{i}:{(\A^1)}^{(p+1)}\times \mathit{Poly}^{n-(p+1)v,r+1} \to {(\A^1)}^{p}\times \mathit{Poly}^{n-pv,r+1}\\
			(a_0,\ldots, a_p), (P_1(z),\ldots, P_r(z)) \mapsto (a_0,\ldots, \hat{a_i}, \ldots, a_{p-1}), \Big((z-a_i)^vP_1(z),\ldots, (z-a_i)^vP_r(z)\Big)
		\end{gather*}That the face maps indeed satisfy the axioms of Definition \ref{def2.10} is explained in Section \ref{sec3}. In particular, $\mathit{Poly}_v^{n,r+1}$ is the space of $\A^1$-indecomposables. The complex points of the space  $\mathit{Poly}_v^{n,r+1}$ (noting that  $\mathit{Poly}_v^{n,r+1}$ is defined over $\mathbb{Z}$ for all $n$), is referred to as $\mathrm{Rat}_n(\P^1,\P^{r-1})$ by Farb-Wolfson, the moduli space of `based holomorphic maps' that take $\infty$ to $[1:\ldots:1]$ (see \cite{Farb2015}). The final result one obtains using Theorem \ref{thm1} is given in Corollary \ref{cor7}.

		\item \emph{The moduli space of degree $n$ maps $C\to \P^r$, $\morn$.} Let $C$ be a smooth projective curve of genus $g\geq 0$ defined over an algebraically closed field $K$ and let $J(C)$ denote the Jacobian of $C$. Let $\pic^n(C)$, which is (noncanonically, by a translation) isomorphic to $J(C)$, denote the space of degree $n$ line bundles on $C$.  A degree $n$ map $C\to \P^r$ is determined by \begin{enumerate}[i.]
			\item a choice of a line bundle $L\in \pic^n(C)$  \item sections $s_0,\cdots, s_r \in H^0(C,L)$ having no common zeroes
		\end{enumerate} whence we have \begin{gather*}
			C\to \P^r\\x\mapsto [s_0(x):\cdots : s_r(x)].
		\end{gather*} Let $\morn$ denote the moduli space of all degree $n$ maps $C\to \P^r$. Define $X_n$ by \begin{gather*}
			X_n:=\{L,[s_0:\ldots:s_r]: L \in \pic^n(C), s_i\in H^0(C,L) \text{ for all }i\}.
		\end{gather*}When $n\geq 2g$ (for $g\geq 2$, even $n\geq 2g-1$ works for our purposes), by the Riemann-Roch theorem $\dim H^0(C,L) =n-g+1$ for all $L\in \pic^n(C)$, which makes $X_n$ the projectivisation of a rank $(r+1)(n-g+1)$ vector bundle on $\pic^n  (C)$ and $\morn\subset X_n$ is Zariski open dense: \footnote{For $n\leq 2g-2$ the description of $X_n$ as the projectivisation of a vector bundle on $\pic^n(C)$ no longer holds; it has been the subject of intense study for decades, under the name of Brill-Noether theory (for a through introduction see \cite[Chapter V]{ACGH}).}

		\begin{center}	\begin{tikzcd}[column sep=huge,row sep=large]
				\P(H^0(C,L)^{r+1}) \cong \P^{(n-g+1)(r+1)-1}\arrow[hookrightarrow]{r}{} & X_n \arrow{d}{}\\
				& \mathit{Pic}^n(C)
		\end{tikzcd} \end{center}
		
		Now $X_n$ has a natural stratification given by the number of common zeroes of an $(r+1)$-tuple of global sections of some degree $n$ line bundle, which in turn shows, by Definition \ref{def2.10} (see Section \ref{sec3} for details) that  $X_n$ admits a symmetric semisimplicial filtration by $C$, with filter gap $e=1$ and the space of $C$-indecomposables $\morn$.
		
	\end{enumerate} 
	
	As we already mentioned, all the examples above are instances of Definition \ref{def2.10}, of spaces admitting symmetric semisimplicial filtration, and the computation of their (stable) cohomology is covered in Section \ref{sec3}. However, the cohomology of $\morn$ is interesting enough that it warrants being recorded in the introduction.
	
	\vspace{3mm}
	
	\noindent \textit{The moduli space of algebraic morphisms of degree $n$ from a genus $g$ curve $C \to \P^r$:}
	Recall that $\morn$ denote the moduli space of degree $n$ morphisms $C\to \P^r$ i.e. \begin{align*}
		\morn:= \Big\{ \big(L, [s_0:\ldots: s_r]\big): \text{ degree of } L= n, \,\,\,\, s_i\in H^0(C,L), \\s_0,\ldots, s_r \text{ have no common zeroes} \Big\}
	\end{align*}
	Furthermore, we denote a vector space spanned by $\{a_1,\ldots, a_k\}$ over $\Q$ by $\Q\{a_1,\ldots, a_k\}$.
	\begin{nmtheorem}\label{thm3}
		Let $C$ be a smooth projective curve of genus $g$, and let $r$ and $n$ be fixed positive integers such that $n\geq 2g$, and $r\leq n-g$ over $\C$. Let $n_0:=n-2g$. Then there exists a second quadrant spectral sequence, which converges to $H^*(\morn; \Qb)$ as an algebra, which has the following description:
		
		\begin{enumerate}
			\item The $E_2$ term is a bigraded algebra. For $p\leq n_0$ the $E_2^{-p,q}$ is the $(p,q)$ graded piece of $$H^*(J(C);\Qb)[h]/h^r \otimes \wedge \Qb\{t\}\otimes \Sym \Qb\{\alpha_1,\ldots, \alpha_{2g}\},$$ where:
			\begin{enumerate}[i]
				\item  $H^i(J(C);\Qb)$ has degree $(0,i)$, $h$ has degree $(0,2)$, $t$ has degree $(-1,2r+2)$ 
				and $\alpha_i$ has degree $(-1,2r+1)$ for all $i$,
				\item modulo elements of degree $(-i,j)$ with $i>n_0$
			\end{enumerate}
			\item  Furthermore this is a spectral sequence of mixed Hodge structures, with $\Qb\{\alpha_1,\ldots, \alpha_{2g}\}$ and $\Qb\{t\}$ each carrying a pure Hodge structure of weight $2(r+1)$, and $h$ is of type $(1,1)$. 
		\end{enumerate}  
	\end{nmtheorem}
	\begin{remark}
		Note that the ground field is assumed to be $\C$ in the theorem above. The proof of the theorem works for all fields except the part where we use some results from Brill-Noether theory (see \cite[Chapter IV, Section 2]{ACGH} on universal divisors). It is widely accepted that most, if not all, of the results from Brill-Noether theory that we use, hold over fields of positive characteristics as well. In turn Theorem \ref{thm3} should hold over any algebraically closed field. However, it seems that there might be some gaps in Brill-Noether theory for positive characteristics in literature, so for the sake of being thorough we stick to $\C$.
	\end{remark}
	\noindent Note that when $C=\P^1$, the Jacobian of $\P^1$ is just a point, and the theorem above gives us the following corollary.
	\begin{cor}\label{cor4}
		Let $r$ and $n$ be positive integers satisfying $r\leq n$. Then $$H^*(\mathrm{Mor}_n(\P^1,\P^r);\Q) \cong \Q[h]/h^r\otimes \wedge\Q\{t\}$$ where $t$ has cohomological degree $2r+1$. Furthermore, over a field $\kappa$, with algebraic closure $\overline{\kappa}$, we have an isomorphism of $Gal (\overline{\kappa}/\kappa)$-representations:
		
		\begin{gather*}
			H^i_{\et}(\mathrm{Mor}_n(\P^1,\P^r);\Qb_{\ell}) = \begin{cases}
				\Qb_{\ell}(-j) & i=2j, 0\leq j\leq r-1\\ \Qb_{\ell}(-(j+1)) & i=2j+1, r\leq j\leq 2r-1\\ 0 & \text{ otherwise.} 
			\end{cases}
		\end{gather*}
		
	\end{cor}
	\begin{remark}
		Strictly speaking Corollary \ref{cor4} follows from Theorem \ref{thm3} only when the ground field is $\C$. We actually prove Corollary \ref{cor4} independently, without using Theorem \ref{thm3}, and the proof works whatever the characteristic of the ground field might be. This is because even though both Corollary \ref{cor4} and Theorem \ref{thm3} are applications of Theorem \ref{thm2}, with the Jacobian of $\P^1$ being just a point, we do not require the power of Brill-Noether theory in all its generality for Corollary \ref{cor4} the way we do for Theorem \ref{thm3}.
	\end{remark}
	One should note here that spaces admitting symmetric semisimplicial filtration arise as a special type of $\Deltainj S$ object. In particular, unlike the previous examples, the following is an instance where the cohomology computation is heavily dependent on the $\Deltainj S$ structure of the moduli space under consideration (see Section \ref{sec4} for details); however, the space itself does not satisfy the conditions of admitting a symmetric semisimplicial filtration by a fixed space.
	\paragraph{More context and one more example: The moduli space of smooth sections of a $\grd$.}
	A \emph{linear series} (or \emph{system}) is a vector subspace of the vector space of global sections of a line bundle on a smooth projective curve. A linear system $V$ on a smooth projective curve $X$ is called a $\grd$ if $V\subset H^0(X,L)$, where $L$ is a degree $d$ line bundle on $X$ and $V$ is a complex $(r+1)$ dimensional vector space. A $\grd$, say $V$, is \emph{$m$-very ample} if for every effective divisor $\xi \in X$ of degree $m+1$, we have that $$\dim\V(-\xi) = r+1-(m+1),$$ where $\V(-\xi):= H^0(X,L\otimes \mathcal{O}(-\xi))\cap \V$.
	Though the following result cannot be obtained as a corollary to Theorem \ref{thm1}, the basic principles of the proof of Theorem \ref{thm1} hold almost verbatim to prove the following on the (stable) cohomology of the moduli space of smooth sections of a $\grd$ that is $m$-very ample. 
	\begin{nmtheorem}\label{thm5}
		Let $X$ be a smooth projective curve of genus $g$ over $\C$. Let $\V$ be a linear system on $X$ of type $\grd$; moreover let $\V$ be $m$-very ample. Define $\V^{\circ} \subset \V$ to be the locus of smooth sections in $\V$. Then for all $i\leq {m-1\over 2}$ the following holds:
		\begin{align*}
			H^i(\V^{\circ};\Qb) \cong \begin{cases}
				\Sym^{p-2}H^1(X;\Qb)(-(p-1)) \oplus \Sym^p H^1(X;\Qb)(-p) & i = 2p\\
				\Sym^{p-1} H^1(X;\Qb) (-(p-1)) \oplus \Sym^p H^1(X;\Qb)(-(p+1)) & i=2p+1.
			\end{cases}
		\end{align*}
	\end{nmtheorem}

	\begin{remark}
		Note that the requirement $r\leq n-g$ is not a major restriction. For a line bundle $L$ over $C$ of degree $n$, satisfying $n\geq 2g$, $\dim H^0(C,L) =n-g+1$ by the Riemann-Roch theorem; so when $r>n-g$, any map $C\to \P^r$ factors through the \emph{complete linear system of $L$} i.e. $$C\hookrightarrow \P(H^0(C,L)^*) \hookrightarrow \P^r$$ where the second embedding is merely a linear embedding. 
	\end{remark}

	\textbf{Relation to other results.} In each of the examples discussed here, the `natural' dense open subsets i.e the `space of indecomposables' are the ones whose topological properties (e.g (co)homology) we are interested in. A lot of work has been done computing the cohomology of such examples. See e.g Church \cite{Church2012}, Totaro (\cite{Totaro}, Farb-Wolfson-Wood \cite{FWW} and the references therein. 
	
	For the second, there are prior results by Segal (\cite{Segal}), Farb-Wolfson (\cite{Farb2015}, for a motivic perspective), Gadish (\cite{Gadish} from the perspective of representation stability) and others. In the same paper (\cite{Segal}), Segal also has results regarding homological stability for Example 3.  Theorem \ref{thm3} and Corollary \ref{cor4} are algebro-geometric and arithmetic generalizations of Segal's \cite{Segal}, on the moduli space of degree $n$ maps $C\to \P^r$ (for the case of $C=\P^1$ over $\C$, see also \cite{GomezGonzales2020} and the references therein). Segal proved that $\morn$ is stably homologous to the moduli space of \emph{continuous maps} $C\to\P^r$ over the ground field $\C$ by a beautiful trick often referred to as `bringing zeroes from infinity', which, of course, works in the analytic topology and does not allow us to keep track of the weights, unlike Theorem \ref{thm3} and Corollary \ref{cor4}.
	
	Our method of proof shares certain similarities with Petersen's work in \cite{Petersen2017} and Tommasi's work in \cite{Tommasi2014} insofar as all have essentially the same root- Deligne's theory of cohomological descent. In particular, \cite[Theorem 3.3]{Petersen2017} computes the cohomology of a simplicial space whose face maps are closed embeddings; and in \cite{Tommasi2014}, Tommasi constructs an augmented semisimplicial space to compute the cohomology of the moduli space of smooth hypersurfaces. 
	
	It has been brought to our attention that Tommasi (\cite{Tommasi}) is also studying the moduli space of smooth sections of a line bundle over a smooth projective curve with the goal of computing some stable cohomology, sans the notion of $m$-very ampleness. A similar topic has been studied by I. Banerjee in \cite{IshanBanerjee} that relates the integral cohomology of the moduli space of sections of a line bundle with certain commutator subgroups of the surface braid group. Parallel to this, Aumonier (\cite{Aumonier2021}) used, like Tommasi, the Vassiliev spectral sequence, and homotopy theoretic methods to show that these moduli spaces are rationally cohomologous (stably) to the moduli space of continuous sections. In this paper, however, we bypass the combinatorial complexities that are involved in Vassiliev's spectral sequence and moreover have the added advantage that our methods are algebraic, and in turn constantly keep track of the weights.
	
	Finally, to the best of my knowledge, the fact that all of these examples can be studied under the same framework joined by a common thread- the property of being a $\Delta S$ or a $\Deltainj S$ object, has not been addressed in the literature, neither has the symmetric simplicial category been exploited to study these examples before. A beautiful paper (with very different goals) that is worth mentioning at this juncture is Dyckherhoff-Kapranov's work on ribbon graphs (see \cite{DK2015})- they use certain crossed simplicial groups to describe the combinatorics that marked surfaces with $G$-structures come equipped with.
	
	Apart from the applications of Theorems \ref{thm1} and \ref{thm2} discussed later in this paper, some immediate consequences of Theorem \ref{thm1} are Theorem A and Corollary B of \cite{Banerjee2021}, which disproves (in Theorem A) Conjecture G' posed by Vakil and Wood in \cite{VakilWood13}, and proves a strengthening of another (Conjecture H' of \cite{VakilWood13}, Corollary B of \cite{Banerjee2021}). The conjectures are centred on certain locally closed subspaces of $\Sym^n(\P^1)$ and the author, in \cite{Banerjee2021}, gives (counter)examples to the  principle of Occam's razor for Hodge structures.
	
	\subsection*{Acknowledgements}
	I am very grateful to Benson Farb for his patient guidance and unconditional support.  His invaluable comments have been instrumental in the improvement of this paper. My deepest thanks to Patrick Brosnan, Lei Chen, Izzet Coskun, Zhiyuan Ding, Quoc Ho, Peter Scholze and Craig Westerland for helpful conversations. I am also thankful to Jesse Wolfson for his thoughtful suggestions, to Ben O'Connor and Alexis Aumonier for their feedback, to Nir Gadish and Claudio Gómez-Gonzáles for their help with the references. Finally a warm thanks to the anonymous referee whose careful reading and some suggested edits made the paper significantly more readable.
	
	\section{The symmetric simplicial group $\fS_{\bullet}$}\label{sec2}
	In this section our goal is to prove Theorems \ref{thm1} and \ref{thm2}. To this end, we first collect some basics on \emph{symmetric (semi)simplicial space}, following \cite{FiedorowiczLoday1991}, then study certain instances of $\Deltainj S$-spaces and locally constructible sheaves, and finally use them to prove theorems \ref{thm1} and \ref{thm2}.
	
	\subsection{\small Generalities on the category $\Delta S$}
	\begin{defn}\label{def2.1}
		Let $[p]:= \{0,\ldots, p\}$ be an ordered set (in the obvious way) with $(p+1)$ elements. Let $\Delta$ be the \textbf{simplicial category} with these objects $[n]$ and morphisms given by monotone maps of ordered sets $[n]\to [m]$. The morphisms of $\Delta$ are generated by the \textbf{face maps} $$f^{\Delta}_j:[p- 1] \to [p]$$ that misses $j$ and the \textbf{degeneracy maps} $$s^{\Delta}_j: [p + 1] \to [p]$$ that hits $j$ twice, $j = 0 \ldots p$. The subcategory $\Deltainj \subset \Delta$ contains all its objects, but only the injective monotone maps.
		The \textbf{symmetric simplicial category}, which we denote by $\Delta S$, is a small category with the following structure:
		\begin{enumerate}[i.]
			\item the objects of $\Delta S$ are $[p]$, $p \geq 0$, 
			\item 	 $\Delta S$ contains $\Delta$ as a subcategory,
			\item $\Aut_{\Delta S}[p]= S^{\mathit{op}}_{p+1}$  (opposite group of $S_{p+1}$),
			\item  any morphism in $\Delta S$ can be uniquely written as a composite $\phi. g$ where $\phi\in \mathrm{Hom}_{\Delta}([p],[m])$ and $g\in S^{\mathit{op}}_{p+1}$.
		\end{enumerate}
		The \textbf{symmetric semisimplicial category} $\Deltainj S\subset \Delta S$ contains, as objects, those in $\Delta S$, and as morphisms, only the injective maps in $\Delta S$. 
	\end{defn}
	
	Define, for each $p$, $$\fS_p:= S_{p+1},$$ the symmetric group on $(p+1)$ letters, and we call $\fS_{\bullet}$ the \emph{symmetric simplicial group}, where $\fS_{\bullet}:=\{\fS_p\}_{p\geq 0}$. In Lemma \ref{lem2.3} and Proposition \ref{prop2.4} we will see that $\fS_{\bullet}$ is not just a collection of symmetric groups; but rather, it is a \emph{simplicial set} whose face and degenracy maps satisfy additional relations. 
	
	Multiplication in $\fS_p$ is the usual group multiplication in $S_{p+1}$; the composition of $g$ and $h$ in  $\Aut_{\Delta S}[p]$ is also given by the group multiplication i.e. $g\circ h= hg$. For all $m,p\geq 0$, and for all $g\in \Aut_{\Delta S}[n]$ and $\phi \in \mathrm{Hom}_{\Delta}([m],[p])$, thanks to the last axiom we have the following commutative diagram: 
	
	\begin{center}
		\begin{tikzcd}[row sep=3.3em, column sep=3.3em]
			{[m]} \arrow[r, "\phi"] \arrow[d, swap, "\phi^*(g)"] \arrow[dr, dashed, "g_{\circ}\phi"]
			& {[p]} \arrow[d, "g"] \\ {[m]} \arrow[r, swap, "g^*(\phi)" ]
			& {[p]} \end{tikzcd}
	\end{center}
	
	\noindent for an unique $\phi^*(g) \in  \Aut_{\Delta S}[m]$ and an unique $g^*(\phi)\in \mathrm{Hom}_{\Delta}([m],[p])$, which defines the composition $g_{\circ}\phi$. Note that $\Delta S$ is naturally equipped with face and degeneracy maps, and we continue to denote them by $f^{\Delta}_j$ and $s^{\Delta}_j$ respectively.
	
	\begin{remark}
		\textup{The symmetric simplicial group is a special case of something much more general- those are called \textup{crossed simplicial groups} (see \cite[Definition 1.1]{FiedorowiczLoday1991}). One such instance arises in the case of Conne's cyclic homology, whose objects are $[p]$, like in $\Delta$, but the automorphism groups are cyclic groups. Other examples of crossed simplicial groups include that formed by the braid groups,  the dihedral groups, the hyperoctahedral groups etc. For a complete treatment of crossed simplicial groups see \cite{FiedorowiczLoday1991}.}
	\end{remark} 
	
	Some important observations on $\fS_{\bullet}$ before we move on to defining objects on the category $\Delta S$. Recall that a \emph{simplicial set} is a functor $T: \Delta^{\mathrm{op}}\to \mathbf{Sets}$, which we often denote by $T_{\bullet}$; i.e. it is a simplicial object in $\mathbf{Sets}$. Unless otherwise stated, for any simplicial object $T_{\bullet}$ in a category $\mathcal{C}$, we will denote its face and degeneracy maps by $f_i^T$ and $s_i^T$ i.e. $$f_i^T:= T(f^{\Delta}_i), \,\,\,\,\, s_i^T:= T(s^{\Delta}_i).$$
	
	\begin{lemma}\label{lem2.3}(\cite[Lemma 1.3]{FiedorowiczLoday1991})
		The symmetric simplicial group	$\fS_{\bullet}$ is a simplicial set given by \begin{align*}
			\fS:  &\Delta^{\mathrm{op}}\to \mathbf{Sets} \\ &[p]\mapsto \fS_p.
		\end{align*}
	\end{lemma}
	\begin{proof} 
		The functor  \begin{align*}
			\fS:  &\Delta^{\mathrm{op}}\to \mathbf{Sets} \\ &[p]\mapsto \fS_p.
		\end{align*} is well-defined by Axiom iv, Definition \ref{def2.1}. Indeed, for $$\phi:[m]\to [p]$$ in $\Delta$, we have a partition of the ordered set $[m]$ as $$[m] = \sqcup_{0\leq i\leq p} \phi^{-1}(i),$$ which in turn defines an unique $$\phi^* :\fS_p\to \fS_m,$$ where $\phi^*(g) \in \fS_m$ permutes elements of $[m]$ by permuting the $p+1$ partition blocks $\phi^{-1}(i)$, $0\leq i\leq p$, by $g\in\fS_p$, while respecting the internal ordering each of the blocks $\phi^{-1}(i)$ posses by virtue being a subset of the ordered set $[m]$.
	\end{proof}
	In the following proposition we record some relations the face maps $f_j^{\fS}$ and degeneracy maps $s^{\fS}_j$ of the simplicial set $\fS_{\bullet}$ satisfy, which would be used later to give a simple characterization of $\fS_{\bullet}$-objects in terms of simplicial objects satisfying additional relations (see Lemma \ref{lem2.6}).
	\begin{prop}\label{prop2.4} For all $p$, and all $\sigma\in \fS_p$, the face maps  $f_j^{\fS}$ and the degeneracy maps $s_j^{\fS}$ of the simplicial set $\fS_{\bullet}$ are such that--
		\begin{enumerate}[i.]
			\item the following relations hold: \begin{gather} \label{eq2.1}
				f_j^{\fS}(\sigma \sigma') = f_j^{\fS}(\sigma) \circ f_{\sigma^{-1} (j)}^{\fS}(\sigma') \\ s_j^{\fS}(\sigma \sigma') = s^{\fS}_j(\sigma) \circ s^{\fS}_{\sigma^{-1}(j)}(\sigma')\nonumber
			\end{gather}
			\item the following diagrams commute:
			\begin{equation}\label{eq2.2}
				\begin{tikzcd}
					{[p-1]} \arrow[r, "f^{\Delta}_{\sigma^{-1}(j)}"] \arrow[d, swap, "f_j^{\fS}(\sigma)"]
					& {[p]} \arrow[d, "\sigma"] \\ {[p-1]} \arrow[r,swap, "f^{\Delta}_j" ]
					& {[p]} \end{tikzcd}
				\hspace{40pt} \begin{tikzcd}
					{[p+1]} \arrow[r, "s^{\Delta}_{\sigma^{-1}(j)}"] \arrow[d, swap, "s_j^{\fS}(\sigma)"]
					& {[p]} \arrow[d, "\sigma"] \\ {[p+1]} \arrow[r, swap, "s^{\Delta}_j" ]
					& {[p]} 
				\end{tikzcd}
			\end{equation}
		\end{enumerate}
	\end{prop}
	\begin{proof}
		The statement of this proposition is a special case of Proposition 1.7 of \cite{FiedorowiczLoday1991}. Plugging in $\fS_{\bullet}$ in place of more general crossed simplicial groups `$G_{*}$' in the proof \cite[Proposition 1.7]{FiedorowiczLoday1991} proves the relations above.
	\end{proof}
	
	Just as the main power of $\Delta$ (and $\Deltainj$) lie in encoding the combinatorial information of objects in a category $\mathcal{C}$ in a succinct fashion by considering functors from $\Delta$, the strength of $\Delta S$ (and, of course, $\Deltainj S$) lie in throwing the extra structure provided by the action of $\{S_p\}_{p\in \mathbb{N}}$ into the mix. 
	\begin{defn}\label{def2.5}
		Let $\mathcal{C}$ be a category. A \textbf{symmetric simplicial object} (or a \textbf{$\fS_{\bullet}$-object}) is a functor $T: (\Delta S)^{\mathrm{op}}\to \mathcal{C}.$ A \textbf{symmetric semisimplicial object} (or \textbf{$\fS_{\bullet, \mathrm{inj}}$-object}) is defined likewise.
	\end{defn}
	We follow the conventional notation from the standard simplicial case: we denote such a functor simply by $T_{\bullet}$. Also, we write $$\fS_{p}\times T_p\to T_p$$ to denote natural action of $S_{p+1}$ on $T_p$ that comes from $T_{\bullet}$ being a $\Delta S$ (respectively $\Deltainj S$) object, and we will denote $(g,t)$ simply by $gt$ for all $g\in \fS_{p}$, $t\in T_p$ and for all $p\geq 0$. In our cases $\mathcal{C}$ will be the category of topological spaces, or the category of schemes (with \etale topology). We will also consider \emph{symmetric cosimplicial objects} (and symmetric cosemisimplicial objects) in the category of $\Q$-vector spaces over a space.
	
	Starting with the following lemma, for the rest of this subsection we try to understand the bridge between objects over $\Delta$ and $\Delta S$. 
	\begin{lemma}\label{lem2.6}
		The notion of a $\fS_{\bullet}$-object in $\mathcal{C}$ is equivalent to the notion of a simplicial object $T_{\bullet}$ in $\mathcal{C}$ with the following additional structure:\begin{enumerate}[i.]
			\item left group actions $\fS_{p}\times T_{p}\to T_{p}$ for all $p\geq 0$,
			\item face relations  $f_j^{T}(\sigma t) = f_j^{\fS}(\sigma)(f^{T}_{\sigma^{-1}(j)} t)$,
			\item 	degeneracy relations $s_j^{T}(\sigma t) = s_j^{\fS}(\sigma)(s^{T}_{\sigma^{-1}(j)} t)$,
		\end{enumerate}
		In fact it suffices to specify the face and degeneracy relations for the generators
		of $\fS_p$. A $\fS_{\bullet}$-map $\phi_{\bullet}: T_{\bullet} \to  T'_{\bullet}$ is the same thing as a simplicial map such that each $\phi_p : T_p \to T'_p$ is $\fS_{p}$-equivariant.
	\end{lemma}
	\begin{proof}
		This is a special case of Lemma 4.2 of \cite{FiedorowiczLoday1991}. The basic idea behind the proof is that the inclusion $\Delta \subset \Delta S$ defines for each $\fS_{\bullet}$-object $T_{\bullet}$ an underlying simplicial object. Paired with Proposition \ref{prop2.4} the statement follows.
	\end{proof}
	
	Now let $X$ be a space and $\mathrm{Sh}_{\Q}(X)$, where $\Q$ is our `coefficient field', be the abelian category of sheaves of $\Q$-modules on $X$.  For sheaves $A, B\in \mathrm{Sh}_{\Q}(X)$ let $\Shom(A, B) \in  \mathrm{Sh}_{\Q}(X)$ denote the internal hom. A \emph{symmetric cosimplicial sheaf}, also called an \emph{$\fS_{\bullet}$-sheaf} or \emph{$\Delta S$-sheaf} is a functor $$\mathcal{F}: \Delta S\to \mathrm{Sh}_{\Q}(X).$$ We will denote the constant sheaf supported on $X$ by $\underline{\Q}_X$. Note that a $\fS_{\bullet}$-sheaf $\mathcal{F}_{\bullet}$ has a natural structure of a cosimplicial sheaf because $\Delta \subset \Delta S$. Also observe that $\mathcal{F}_n$ is naturally a sheaf of $\underline{\Q}_{X}[\fS_{p}]$-modules. 
	Let $\big({\mathcal{F}}_{p}\otimes  \sgn_{p+1}\big)^{S_{p+1}}$ denote symmetric group invariants of the sheaf $\mathcal{F}_p$ under the permutation action of $S_{p+1}$ twisted by the sign. 
	
	\begin{lemma}\label{lem2.7}
		The sheaf $\Sext^{\,\,*}_{\Delta S}(\underline{\Q}_X, \mathcal{F}_{\bullet}) $ is isomorphic to the homology of the complex of coinvariants $$\big({\mathcal{F}}_{p}\otimes  \sgn_{p+1}\big)^{S_{p+1}}, d$$ where $d= \sum (-1)^i f_i$. 
	\end{lemma}
	\begin{proof}
		The statement is a special case of Corollary 6.10 of \cite{FiedorowiczLoday1991}. The idea is roughly as follows. In \cite{FiedorowiczLoday1991} they choose a certain biresolution of $\underline{\Q}_X$  which gives rise to a bicomplex computing $\Sext^{\,\,*}_{\Delta S}(\underline{\Q}_X, \mathcal{F}_{\bullet})$. The naive filtration by row gives a spectral sequence with the cohomology of $S_{p+1}$ with coefficients in $\big({\mathcal{F}}_{p}\otimes  \sgn_{p+1}\big)^{S_{p+1}}$ on the $E_1$ page and the differentials are naturally given by the alternating sum of the face maps. And this spectral sequence converges to $\Sext^{\,\,*}_{\Delta S}(\underline{\Q}_X, \mathcal{F}_{\bullet})$. 
	\end{proof}
	\noindent Recall that since $\mathcal{F}_{\bullet}$ is a cosimplicial sheaf, the cohomology of the complex $$\mathcal{F}_0\to \mathcal{F}_1\to \cdots$$ with differentials given by the alternating sum of the face maps $d=\sum (-1)^i f_i$ is given by $\Sext^{\,\,*}_{\Delta}(\underline{\Q}_X,\mathcal{F}_{\bullet})$. 
	\begin{lemma}\label{lem2.8}
		If $\mathcal{F}_{\bullet}$ is an $\fS_{\bullet}$-sheaf, then the canonical map $$\Sext^n_{\Delta }(\underline{\Q}_X, \mathcal{F}_{\bullet}) \to \Sext^n_{\Delta S}(\underline{\Q}_X, \mathcal{F}_{\bullet}) $$ is an isomorphism for all $n$.
	\end{lemma}
	\begin{proof}
		For a proof see Theorem 6.16 of \cite{FiedorowiczLoday1991}.
	\end{proof}
	\noindent The main takeaway from this subsection is the following proposition:
	
	\begin{prop}\label{prop2.9}
		Let $\mathcal{F}$ be a $\fS_{\bullet}$-sheaf. The surjection of complexes $$\Big(\mathcal{F}_p, d\Big) \to \Big(\big({\mathcal{F}}_{p}\otimes  \sgn_{p+1}\big)^{S_{p+1}}, d\Big),$$ where $d=\sum (-1)^i d_i$,  is a quasi-isomorphism.
	\end{prop}
	\begin{proof}
		Follows immediately from lemmas \ref{lem2.7} and \ref{lem2.8}. Indeed, by Lemma \ref{lem2.8} we have $$\Sext^p_{\Delta }(\underline{\Q}_X, \mathcal{F}_{\bullet}) \xrightarrow{\cong} \Sext^p_{\Delta S}(\underline{\Q}_X, \mathcal{F}_{\bullet}) $$ and by Lemma \ref{lem2.7} we have $$\Sext^p_{\Delta S}(\underline{\Q}_X, \mathcal{F}_{\bullet}) \xrightarrow{\cong} \Sext^p(\underline{\Q}_X, \big({\mathcal{F}}_{p}\otimes  \sgn_{p+1}\big)^{S_{p+1}}).$$
	\end{proof}
	\subsection{\small Spaces admitting symmetric (semi)simplicial filtration}
	A particular subclass of $\Deltainj S$-objects deserves special attention simply because of the sheer number of examples that fit into it (see the examples in the introduction). 
	Let $M$ and $\{X_n\}_{n\in \mathbb{N}}$ be spaces, $X:=\sqcup X_n$, and $e$ a positive integer. Recall that ${S}_p$ denotes the symmetric group on $p$ elements. For an element $\sigma \in S_p$ we let $\sigma$ denote the automorphism on $M^p$ induced by permuting the factors as well. 
	For any space $X$ let $\mathit{id}_X$ denote the identity map on $X$. \begin{defn}\label{def2.10}
		We say that $\{X_n\}_{n\in\mathbb{N}}$  \textbf{admits a symmetric semisimplicial filtration by $M$} with \textbf{filter gap} $e>0$, if for all $0\leq i\leq p\leq \frac{n}{e}$ there are proper finite morphisms called the \textbf{face maps},  with the $i^{th}$ face map defined as  \begin{align*}
			f_{i}:M^{p+1} \times X_{n-e(p+1)} \to M^{p}\times X_{n-ep},
		\end{align*}  and satisfying the following axioms. 
		\begin{enumerate}[i.]
			\item (\textbf{Semisimplicial identity}) For all $i<j$ the following hold: \begin{align}\label{eq2.3}
				f_i \circ f_{j} = f_{j-1} \circ f_{i} 
			\end{align}
			\item  (\textbf{Symmetric condition}) Given $\sigma \in S_{p+1}$ and $x\in X_{n-e(p+1)}$ and for each $i\geq 0$, there exists a unique $d_i(\sigma) \in S_{p}$ and $x'\in X_{n-ep}$ such that \begin{equation*}
				f_i\big(\sigma(t_0, \cdots, t_p),x\big) = d_i(\sigma) f_{\sigma^{-1}(i)}\big((t_0,\cdots , t_p),x'\big) 
			\end{equation*} i.e. the following diagram commutes: 
			
			\begin{equation}\label{eq2.4}
				\begin{tikzcd}[row sep=3.5 em]
					{M^{p+1}\times X_{n-e(p+1)}} \arrow[r, "f_i"] \arrow[d, "\sigma\times id_{X_{n-e(p+1)}}"]
					& {M^p\times X_{n-ep}} \arrow[d, "d_i(\sigma)"] \\ {M^{p+1}\times X_{n-e(p+1)}} \arrow[r, "f_{\sigma^{-1}(i)}" ]
					& {M^p\times X_{n-ep}} \end{tikzcd}
			\end{equation}

			\item  (\textbf{Equalizer condition}) Let
			$$\pi_p: M^{p+1} \times X_{n-e(p+1)} \to X_n$$
			be defined by  \begin{align*} \pi_p: = f_{0}\circ f_{1}\circ \cdots \circ f_{p-1} \circ f_{p}.\end{align*} 
			If $(z_0,x_0),\ldots, (z_{p},x_{p}) \in M\times X_{n-e}$ are such that $f_{0}(z_i,x_i) = f_{0}(z_j,x_j) = x$ for some $x\in X_n$ and  all $0\leq i\leq j\leq p$, then there exists a unique $y\in X_{n-e(p+1)}$ such that \begin{align}\label{eq2.5}
				\pi_p((z_0,\ldots, z_{p}), y) = x.
			\end{align} \item   (\textbf{Embedding condition}) For all $z\in M$ the morphisms \begin{align}\label{eq2.6}
				f_{0}(z,\underline{\hspace{3mm}}):X_{n-e} \to X_n
			\end{align} are embeddings (in the relevant category).			
		\end{enumerate}
		We call $$U_n:= X_n - f_0(M\times X_{n-e})$$  the \textbf{space of $M$-indecomposables}. 	\hfill $\square$
	\end{defn}  
	\noindent Recall from the introduction that $X=\sqcup_n X_n$. Now observe that the first two axioms imply $X$ a module over the graded commutative monoid $\Sym M$, where $M$ has grading $e$ and $X$ is naturally graded by $n$. The first two conditions also naturally make $T_{\bullet}$, defined by $$T_{p} :=  M^{p+1}\times X,$$ a $\Deltainj S$ space. Note that for each $p$ the space $T_p$ is naturally graded by $n$, and for each $n$ we can define $T_{\bullet, n}$ to be the $n^{\mathrm{th}}$ graded piece of $T_{\bullet}$ i.e. for each $n$ we have: \begin{align*}
		T_{\bullet,n}: & \Deltainj S\to \mathbf{Spaces}\\ & [p]\mapsto M^{p+1} \times X_{n-e(p+1)}.
	\end{align*}
	
	\noindent  Indeed, the commutative diagram \eqref{eq2.4} is precisely a restatement of the face relations in Lemma \ref{lem2.7} (compare it with commutative diagram \eqref{eq2.2}) once we put the face maps $f_i$ of Definition \ref{def2.10} in place of $f_i^T$ in Lemma \ref{lem2.6}.

	Let $\F$ be a $\Q$-sheaf on $X$ and let $\F_n$ denote its restriction to $X_n$.
	Observe that for all $p\geq 0$ and each $n$, the sheaves $\pi_{p_*}\pi_p^*\F_n$ are equipped with the permutation action of ${S}_{p+1}$. It is easy to see that $T_{\bullet,n}$ being an $\fS_{\bullet,\mathrm{inj}}$-space for each $n$ implies that  $\pi_{p_*}\pi_p^*\F$ is an $\fS_{\bullet,\mathrm{inj}}$-sheaf. As before, let $\pi_{p_*}\pi_p^*\F \otimes \sgn_{p+1}$ denote the sheaf $\pi_{p_*}\pi_p^*\F$ with the permutation action of $S_{p+1}$ twisted by a sign. Taking its coinvariants under $S_{p+1}$ we get a complex given by $$C^p(\F): =  (\pi_{p_*}\pi_p^*\F\otimes \sgn_{p+1})^{{S}_{p+1}}.$$  The face maps $f_{i}$ in Definition \ref{def2.10} induce face maps on the cosimplicial sheaf $\pi_{{\bullet}_*}\pi_{\bullet}^*\F$, which we denote by $f_i^*$. The differentials in this complex are given by the alternating sum of the face maps $f_i^*$ on $\pi_{\bullet_*}\pi^*_{\bullet}\F$, like we had for Lemma \ref{lem2.6}, i.e. $$d:C^{p}(\F) \to C^{p+1}(\F)$$ by $$d:= \sum_i(-1)^i f_i^*$$ resulting in a complex of sheaves on $X_n$ given by $C^{\bullet}(\F)$. Note that we can define the complex $C^{\bullet}(\F_n)$ exactly the same way and  $$C^{\bullet}(\F) = \bigoplus_n C^{\bullet}(\F_n).$$
	
	Our next lemma is the sheaf-theoretic analogue of Theorem \ref{thm1} and is the key step towards proving it.  
	
	\begin{lemma}\label{lem2.11}
		Let $X=\sqcup_n X_n$ admit a symmetric semisimplicial filtration by $M$. Let $e>0$ be the filter gap and $U=\sqcup_n U_n$ be the space of $M$-indecomposables. For each $n\in \N$ the complex $(C^{\bullet}(\F),d)$, which reads as \begin{align}
			0\to j_{!}j^*\F\to \F\to \pi_{0_*}\pi_0^*\F \to (\pi_{1_*}\pi_1^*\F\otimes \sgn_{2})^{{S}_{2}} \cdots \to\nonumber \\ \cdots \to  (\pi_{p_*}\pi_p^*\F\otimes \sgn_{p+1})^{{S}_{p+1}}\to \cdots \label{eq2.7}
		\end{align} is exact.
		
		In particular, plugging in $\F=\Q_X$ and restricting to $X_n$, the complex $(C^{\bullet}(\Q_{X_n}),d)$ is the following exact sequence of sheaves: \begin{align}\label{eq2.8}
			0\to j_{!}j^*\underline{\Q}_{U_n} \to \underline{\Q}_{X_n} \to \pi_{0_*}\underline{\Q}_{T_0} \to (\pi_{1_*}\underline{\Q}_{T_1}\otimes \sgn_{2})^{{S}_{2}} \cdots \to \nonumber \\ \cdots \to  (\pi_{p_*}\underline{\Q}_{T_p}\otimes \sgn_{p+1})^{{S}_{p+1}}\to \cdots 
		\end{align} 
	\end{lemma} 
	
	\begin{proof}
		We give two proofs of \eqref{eq2.7}: one using the fact that $C^{\bullet}(\F_n)$ is a $\Deltainj S^{\mathrm{op}}$-sheaf, and the other, in the case $\F = \underline{\Q}_{X}$ (which is the case we would actually need to prove our theorems) involves essentially checking by hand that \eqref{eq2.8} is exact at the level of stalks.
		
		\begin{enumerate}[\emph{Method} 1:]
			\item Note that for each $n$, the locally constant sheaf $\pi_{\bullet_*}\pi_{\bullet}^*\F_n$ is a $\Deltainj^{\mathrm{op}}$-sheaf 
			$$C^p(\F_n)= (\pi_{p_*}\pi_p^*\F_n\otimes \sgn_{p+1})^{{S}_{p+1}}$$ for all $p\geq 0$. Let $$F: \mathrm{Fun}(\Delta S, \mathrm{Sh}_{\Q}(X_n))\to \mathrm{Fun}(\Deltainj S, \mathrm{Sh}_{\Q}(X_n))$$ denote the forgetful functor taking symmetric cosimplicial objects to symmetric cosemisimplicial ones forgetting degeneracies, and let $$F': \mathrm{Fun}(\Deltainj S, \mathrm{Sh}_{\Q}(X_n))\to \mathrm{Fun}(\Delta S, \mathrm{Sh}_{\Q}(X_n))$$  denote its left adjoint, called freely adding degeneracies. We abuse notation and use $F$ and $F'$ to denote similar functors for $\Delta^{\mathrm{op}}$ and $\Deltainj^{\mathrm{op}}$ sheaves as well.
			The $\Deltainj S^{\mathrm{op}}$ sheaf $\pi_{\bullet_*}\pi^*_{\bullet}\F$ gives us a $\Delta S^{\mathrm{op}}$ sheaf $F'\Big(\pi_{\bullet_*}\pi^*_{\bullet} \F\Big)$ by `freely adding degeneracies'.  In turn, the cohomology $H^n\big(C^{\bullet}(\F)\big)$ is isomorphic to the cohomology of the complex whose terms are given by $$\Big( F'\big(\pi_{p_*}\pi^*_{p} \F\big) \otimes \sgn_{p+1}\Big)^{S_{p+1}},$$
			because freely adding degeneracies do not change the cohomology (say, by Dold-Kan). By Lemma \ref{lem2.8}, this is isomorphic to $$\Sext^n_{\Delta}\Big(\underline{\Q}_{X_n}, F'\big(\pi_{\bullet_*}\pi^*_{\bullet} \F\big)\Big)$$ where we consider the underlying simplicial structure of the $\Delta S^{\mathrm{op}}$-sheaf $F'\big(\pi_{\bullet_*}\pi^*_{\bullet} \F\big)$, and that in turn is isomorphic to cohomology of the the corresponding $\Deltainj^{\mathrm{op}}$-sheaf $$\Sext^n_{\Delta}\big(\underline{\Q}_{X_n}, \pi_{\bullet_*}\pi^*_{\bullet} \F \big),$$ and they vanish, because for each $n$ $$T_{\bullet,n}\to X_n-U_n$$ is proper (in fact finite, because the face maps are finite) and surjective, and thus admits cohomological descent (see \cite[5.3.5(II)]{Deligne1975}), which completes the proof.
			
			\item   If $x\in U_n$ then the stalks of $ j_{!}j^* \underline{\Q}_{X_n}= j_!\underline{\Q}_{U_n}$, and $\underline{\Q}_{X_n}$, both are one dimensional $\Q$-vector spaces each, and the sheaves $\pi_{p_*}\underline{\Q}_{T_p}$  have stalk $0$ at $x$ for all $p\geq 0$. 
			
			So for the rest of the proof we fix $x\in X_n-U_n$. Let $$\pi_0^{-1}(x)= \{(m_0,\ast), \cdots, (m_r,\ast)\}\subset M\times X_{n-e}$$ where $\ast$ denotes not-necessarily-equal elements of $X_{n-e}$, which will not keep track of. Then, first note that $$ \big(j_{!}j^* \underline{\Q}_{X_n}\big)_x= \big(j_!\underline{\Q}_{U_n}\big)_x=0$$ and $$\big(\pi_{p_*}\underline{\Q}_{T_p}\big)_x=0$$ for all $p>r$. So we only have to compute $\big(\pi_{p_*}\underline{\Q}_{T_p}\big)_x$ for $0\leq p\leq r$. Now, the stalk of $\pi_{0_*}\underline{\Q}_{T_0}$ at $x$ is a $\Q$-vector space spanned by a set indexed by $\pi_0^{-1}(x)$, so $$\big(\pi_{0_*}\underline{\Q}_{T_0}\big)_x \cong \Q\{m_0,\ldots,m_r\}.$$
			
			Likewise, $\big(\pi_{1_*}\underline{\Q}_{T_1}\big)_x$ is a $\Q$-vector space spanned by the set $$f_0^{-1}(\pi_0^{-1}(x)) \sqcup f_1^{-1}(\pi_0^{-1}(x)) $$ where $$f_0,f_1: M^2\times X_{n-2e}\to M\times X_{n-e}$$ are the face maps; note that the following sets have the same elements: $$f_1^{-1}(\pi_0^{-1}(x)) = f_0^{-1}(\pi_0^{-1}(x))= \{(m_i,m_j,\ast)\}_{0\leq i,j\leq r}$$ and these sets are nonempty by \eqref{eq2.5} in the Equalizer condition of Definition \ref{def2.10}. For convenience we write tuples of elements of $M$ `multiplicatively', as elements of the tensor algebra $TM:= \sqcup M^n$, which is a associative graded monoid. By this notation, $$\big(\pi_{1_*}\underline{\Q}_{T_1}\big)_x\cong \Q{\{m_im_j: 0\leq i,j\leq r\}}.$$
			The symmetric group $S_2$ acts by swapping the $m_i$ and $m_j$, therefore $$\big(\pi_{1_*}\underline{\Q}_{T_1}\otimes \sgn_{2}\big)^{S_2}_x \cong \Q\{(m_i,m_j)\}_{0\leq i<j\leq r} \cong  \Q\{m_i\wedge m_j:0\leq i<j\leq r\}.$$
			
			Continuing this way, keeping track of the preimages of $x$ under $\pi_0$ followed by various face maps, we get that for all $p\leq r$, the vector space $\big(\pi_{p_*}\underline{\Q}_{T_p}\big)_x$ is isomorphic to a $\Q$-vector space spanned by $$ \{ m_{i_0}\cdot\ldots\cdot m_{i_p})\}_{0\leq i_0,\ldots,i_p\leq r},$$ and therefore \begin{align*}
				\big(\pi_{p_*}\underline{\Q}_{T_p}\otimes \sgn_{p+1}\big)^{S_{p+1}}_x\cong \Q{\{(m_{i_0}\wedge\ldots\wedge m_{i_p}): 0\leq i_0<\ldots<i_p\leq r}\}
			\end{align*}
			
			Having computed the stalks of the terms in \eqref{eq2.8}, we now compute the differential, which are given by alternating sum of the face maps. The differential $$ \big(\underline{\Q}_{X_n}\big)_x \to \big(\pi_{0_*}\underline{\Q}_{T_0}\big)_x $$ is given by $$f_0^*(x) =\sum m_i.$$The differential $$\big(\pi_{0_*}\underline{\Q}_{T_0}\big)_x \to \big(\pi_{1_*}\underline{\Q}_{T_1}\otimes \sgn_{2}\big)^{S_{2}}_x$$ is given by $f_0^*-f_1^*$, where $f_0, f_1$ are the two face maps from $T_1\to T_0$. Noting that $$f_0^*(m_i) = \sum_j m_j\wedge m_i$$ and $$f_1^*(m_i) =m_i\wedge  \sum_j m_j$$ we see that \begin{align*}
				f_0^*-f_1^*: &\big(\pi_{0_*}\underline{\Q}_{T_0}\big)_x \to \big(\pi_{1_*}\underline{\Q}_{T_1}\otimes \sgn_{2}\big)^{S_{2}}_x \\
				& m_i\mapsto 2(\sum_jm_j)\wedge m_i.
			\end{align*} \noindent The same proof applies to show that the differentials are all $$\ast \wedge \sum m_i.$$ Therefore, at the level of stalks \eqref{eq2.8} reads as:
			\begin{align*}
				0\to \Q\xrightarrow{\wedge (\sum m_i)}\Q\{m_0,\ldots, m_r\}\xrightarrow{\wedge (\sum m_i)} \Q\{(m_i\wedge m_j)\}_{0\leq i<j\leq r}\} \xrightarrow{\wedge (\sum m_i)} \cdots\\ \cdots\xrightarrow{\wedge (\sum m_i)}  \Q{\{(m_{i_0}\wedge\ldots\wedge m_{i_p}): 0\leq i_0<\ldots<i_p\leq r}\}\to \cdots\\ \cdots \xrightarrow{\wedge (\sum m_i)} \Q\{m_0\wedge\cdots\wedge m_r\} \xrightarrow{\wedge (\sum m_i)} 0
			\end{align*} which is a Koszul complex and thus exact.
			
		\end{enumerate}	
	\end{proof}

	\begin{remark}
		The reason for giving a proof of Lemma \ref{lem2.11} that exploits only the $\Deltainj S$ structure without appealing to the module structure of $X$ over $\Sym M$ is that the lemma holds even when a $\Deltainj S$ or $\Delta S$ object does not satisfy all the axioms of Definition \ref{def2.10}- evidence at hand is the proof of Theorem \ref{thm3} in Section \ref{sec4}. The space of (global algebraic/holomorphic) sections of a $\grd$ on a smooth projective curve does not satisfy the axioms of Definition \ref{def2.10}; and yet, there's a natural $\Delta S$-space that is homotopy equivalent (in the appropriate category) to the locus of sections with singularities, and a version of Lemma \ref{lem2.11} is the key step towards computing the cohomology of the moduli space of smooth sections of a $\grd$. 
	\end{remark}
	Recall our notations and conventions: $\mathbf{Q}$ denotes $\mathbb{Q}$, the field of rational numbers, or $\mathbb{Q}_{\ell}$, the field of $\ell$-adic numbers as the situation dictates. For any space $X$ (recall that by space we mean locally compact Hausdorff topological space or a quasi-projective algebraic variety over some field), we let $H^*_c(X;\mathbf{Q})$ denote the \etale cohomology with proper supports with coefficients in $\Q_{\ell}$ if $X$ is a quasi-projective algebraic variety, or singular cohomology with compact support with $\Q$ coefficients if $X$ is a topological space. Now, we prove theorems \ref{thm1} and \ref{thm2}.
	\begin{proof}[Proof of Theorem \ref{thm1}]
		Plugging in $\F_n= \underline{\mathbf{Q}}_{X_n}$, we know from Lemma \ref{lem2.11} that $C^{\bullet}(\underline{\Q}_{X_n})$ is a resolution of $j_{!}\underline{\Q}_{U_n}$. Taking cohomology with compact supports we obtain a spectral sequence which reads as $$E_1^{p,q} = H^q_c\Big(X_n, \big({\pi_{p-1}}_*\Q_{T_{p-1}}\times \sgn_{p}\big)^{{S}_{p}}\Big) = \big(H^q_c(M^{p}\times X_{n-ep};\Q)\otimes \sgn_{p}\big)^{{S}_{p}} $$ where we define $T_{-1}:= X_n$ and $S_0$ to be the trivial group. \footnote{We deviate from the standard convention here: for an augmented (semi)simplicial space $T_{\bullet}\to T$, the notation $T_{-1}$ denotes the space $T$. So if we followed the standard convention $T_{-1}$, in our case, should have been $X_n-U_n$, but instead, for convenience, we deviate from what's standard and define $T_{-1} = X_n$.} Applying the Kunneth formula gives: $$E_1^{p,q} =\bigoplus_{l+m=q} \big(H^l_c(M^{p};\Q)\otimes \sgn_{p}\big)^{{S}_{p}}\otimes H^m_c(X_{n-ep};\Q).$$
		Note (for example from \cite{Macdonald1962}) that for $\alpha, \beta \in H^*_c(M;\Q)$ we have $$\alpha\beta - (-1)^{\deg (\alpha)\, \deg(\beta) +1}\beta\alpha =0$$ when taking invariants under the alternating action of ${S}_p$ i.e. \begin{gather*}
			\Big(H^*_c(M;\Q)^{\otimes p} \otimes \sgn_p\Big)^{{S}_p} = H^*_c(M;\Q)^{\otimes p} / \{\alpha\beta - (-1)^{\deg (\alpha)\, \deg(\beta) +1}\beta\alpha =0\}.
		\end{gather*} Therefore, $$E_1^{p,q} = \bigoplus_{l+m=q}\,\,\bigoplus_{i+j=p}  \big(\Sym^i H_c^{\textrm{odd}}(M;\mathbf{Q}) \otimes \mathbb{\wedge}^j H_c^{\textrm{even}}(M;\mathbf{Q})\big)^{(l)} \otimes H^m_c(X_{n-ep};\mathbf{Q}) .$$ In the algebraic setting, when all spaces are quasi-projective varieties over a field $K$, since the face maps are all algebraic morphisms, this spectral sequence is that of $Gal(\overline{K}/K)$ representations.
	\end{proof} 
	
	
	\begin{proof}[Proof of Theorem \ref{thm2}.]
		As before we plug in $\F_n= \underline{\mathbf{Q}}_{X_n}$ to obtain the complex $$j_{!}\underline{\Q}_{U_n}\hookrightarrow C^{\bullet}(\underline{\Q}_{X_n})$$ and it is exact in the category $\mathrm{Sh}_{\Q}(X_n)$ thanks to Lemma \ref{lem2.11}. Let $$\tau: X_n\to \mathrm{pt}$$ be the structure map to a point and let $$\tau_p:=\pi_p\circ \tau: T_p\to \mathrm{pt}$$ denote the respective structure maps for all $p$. We take the global Verdier dual and focus on the resulting complex $$\mathrm{RHom}\big( C^{\bullet}(\underline{\Q}_{X_n}), \underline{\Q}_{X_n}\big).$$ With the naive filtration on the columns of this complex we get a second quadrant $E_1$-page spectral sequence that reads as: \begin{align}
			E_1^{-p,q} = \mathrm{Ext}^q\Big((\pi_{{p-1}_*}\underline{\Q}_{T_{p-1}}\otimes \sgn_{p})^{{S}_{p}}, \underline{\Q}_{X_n}\Big) \implies H^{q+p}(U_n;\Q).
		\end{align}
		Define $N:= \dim_{\mathbb{R}} (X_n)$ if $X_n$ is a smooth orientable manifold, or $N:= 2 \dim_{K} (X_n)$ if $X_n$ is a smooth projective variety over an algebraically closed field $K$. Recall , from the statement of Theorem \ref{thm2}, that $$c(n,p):= 2\big(\dim_{K}X_n - \dim_{K} T_{p-1}\big).$$  Now we compute the $E_1^{-p,q}$ terms by going through the following sequence of steps:
		\begin{flalign*}
			\mathrm{Ext}^q(\pi_{{p-1}_*}\underline{\mathbf{Q}}_{T_{p-1}},\underline{\mathbf{Q}}_{X_n}) &= \mathrm{Ext}^{q-N}(\pi_{{p-1}_*}\underline{\mathbf{Q}}_{T_{p-1}},\underline{\mathbf{Q}}_{X_n}[N]) &&\text{adjusting shifts},\\
			&= \mathrm{Ext}^{q-N}(\pi_{{p-1}_!}\underline{\mathbf{Q}}_{T_p},\underline{\mathbf{Q}}_{X_n}[N])&& \text{$\pi_{p-1}$ finite, $\pi_{{p-1}_*} = \pi_{{p-1}_!}$,}\\
			& = \mathrm{Ext}^{q-N}(\underline{\mathbf{Q}}_{T_{p-1}},\pi_{p-1}^{!}\underline{\mathbf{Q}}_{X_n}[N]) &&\text{ $(\pi_{{p-1}_!}, \pi_{p-1}^!)$ adjoint pair,}\\
			&	=  \mathrm{Ext}^{q-N} (\underline{\mathbf{Q}}_{T_{p-1}},  \pi_{p-1}^! \tau^! \underline{\mathbf{Q}}_{pt}) && \text{$X_n$ smooth, $\Q_{X_n}[N]= \tau^!\underline{\Q}_{pt}$,}\\
			&=  \mathrm{Ext}^{q-N} (\underline{\mathbf{Q}}_{T_{p-1}},  {\tau}_{p-1}^! \underline{\mathbf{Q}}_{pt}) && \tau_{p-1} = \tau\circ \pi_{p-1},\\
			&	= \mathrm{Ext}^{q-N} (\underline{\mathbf{Q}}_{T_{p-1}}, \underline{\mathbf{Q}}_{T_{p-1}}[N-c(n,p)])  && \text{$T_{p-1}$ smooth and} \\ & \text{} && \underline{\Q}_{T_{p-1}}[N-c(n,p)]= \tau_{p-1}^!\underline{\Q}_{pt},\\ 
			&	= \mathrm{Ext}^{q-c(n,p)}(\underline{\mathbf{Q}}_{T_{p-1}}, \underline{\mathbf{Q}}_{T_{p-1}}) && \text{adjusting shifts,}\\&= H^{q-c(n,p)}(T_{p-1}, \mathbf{Q}).
		\end{flalign*}
	\end{proof}
	
	\section{Computing cohomology of some spaces admitting symmetric semisimplicial filtration}\label{sec3}
	In this section we study, in further detail, the impact of Theorem \ref{thm1} on the examples of semisimplicially filtered spaces presented in the introduction. We also prove Theorem \ref{thm5} as an instance of a $\Delta S$ object. All the examples we see in this section are well-known and well-studied by various other methods. However, our Verdier-Deligne inspired approach unifies all of these examples under one framework (which is the property of admitting a symmetric semisimplicial filtration). This, to the best of my knowledge, is new. 
	
	\subsection{\small Unordered configuration spaces}\label{sub3.1} We elaborate on Example 1 from the introduction.
	Recall that $X$ is space, and we define a family of spaces $X_n :=\Sym^nX$ for all $n\in\N$. Also recall the face maps from \eqref{eq1.3}: \begin{gather}
		f_{i}: X^{p+1}\times \Sym^{n-2(p+1)}X \to X^{p} \times \Sym^{n-2p}X \nonumber \\
		\Big( (a_0,\ldots, a_p), \{b_1,\ldots, b_{n-2(p+1)} \}\Big)  \mapsto  \Big((a_0,\ldots, \hat{a_i}, \ldots, a_p), \{ a_i, a_i , b_1,\ldots, b_{n-2(p+1)}\}\Big) \nonumber
	\end{gather} where $\hat{a_i}$ means $a_i$, the $(i+1)^{th}$ factor, removed.
	That the face maps satisfy all the axioms from Definition \ref{def2.10}, with $\uconf_n(X)$ as the space of $X$-indecomposables and $e=2$ as the filter gap, is almost immediate. Therefore, plugging in $M=X$ and $e=2$ in Theorem \ref{thm1}, the spectral sequence from \eqref{eq1.1} reads as:
	\begin{cor}\label{cor6}
		Let $X$ be a locally compact Hausdorff topological space. Then there exists a spectral sequence \begin{gather}
			E_1^{p,q} = \bigoplus_{i+j=p}\bigoplus_{l+m=q}\Big(\Sym^iH^{\textrm{odd}}_c(X;\mathbb{Q})\otimes \Lambda^j H^{\textrm{even}}_c(X;\mathbb{Q})\Big)^{(l)}\otimes H^m_c(\Sym^{n-2p}X;\mathbb{Q})\nonumber \\\implies H^{p+q}_c(\uconf_n X;\mathbb{Q}).\label{eq3.1}
		\end{gather}	where the differentials are given by alternating sum of the pullbacks on cohomology induced by the face maps: \begin{gather*}
			d_1^{p,q}: E_1^{p,q} \to E_1^{p+1,q} \\ d_1^{p,q}:= \sum_{i=0}^{p-1}(-1)^i f_i^*.
		\end{gather*}with differentials given by the alternating sum of the pullbacks on cohomology with compact supports induced by the face maps: \begin{gather*}
			\big((\alpha_1 \cdots \alpha_i )\otimes (\beta_1\wedge \cdots \wedge \beta_j)\big) \otimes  \big((\beta'_1 \cdots \beta'_{j'} )\otimes (\alpha'_1\wedge \cdots \wedge \alpha'_{i'})\big) \mapsto \\ \sum_{1\leq r<s\leq i'} (-1)^{r+s}\big((\alpha_1 \cdots \alpha_i (\alpha'_{r}+\alpha'_s))\otimes (\beta_1\wedge \cdots \wedge \beta_j)\big) \otimes  \big((\beta'_1 \cdots \beta'_{j'} )\otimes (\alpha'_1\wedge\cdots \wedge \widehat{\alpha'_r}\cdots\wedge \widehat{\alpha'_s} \wedge\cdots \alpha'_{i'})\big) 
		\end{gather*}
		where $i+j=p$, $i'+j' = n-2p$ and $$\alpha_1, \ldots, \alpha_i, \alpha'_1,\ldots, \alpha'_{i'} \in H_c^{{\textrm{odd}}}(X)$$ and $$\beta_1, \ldots, \beta_j, \beta'_1,\ldots, \beta'_{j'} \in H_c^{\textrm{even}}(X).$$
	\end{cor}
	
	In the particular case when $X=\C$, the complex numbers, we think of the family $X_n$ can be interpreted as the space of all monic polynomials of degree $n$ over $\C$, which we denote by $(\C[x])_n$. The face maps can be rewritten in terms of multiplication of polynomials: \begin{gather}
		f_{i}: \C^{p+1}\times (\C[x])_{n-2(p+1)}\to \C^{p} \times (\C[x])_{n-2p}\nonumber \\ \Big( (a_0,\ldots, a_p), P(x) \Big)\mapsto  \Big((a_0,\ldots, \hat{a_i}, \ldots, a_p), (x-a_i)^2P(x)\Big)
	\end{gather} As before, $\uconf_n(\C)$, the subspace of square-free polynomials is the space of $\C$-indecomposables. Then Corollary \ref{cor6} \eqref{eq3.1}gives us: \begin{gather} 
		E_1^{p,q} =\bigoplus_{l+m=q}\,\,\bigoplus_{i+j=p} \big(\Sym^i H_c^{\textrm{odd}}(\C;\mathbb{Q}) \otimes \Lambda^j H_c^{\textrm{even}}(\C;\mathbb{Q})\big)^{(l)} \otimes H^m_c(\Sym^{n-2p}\C;\mathbb{Q}) \nonumber \\ \implies H^{p+q}_c(\uconf_n(\C);\mathbb{Q})
	\end{gather} Noting that $H^2_c(\C;\mathbb{Q}) = \mathbb{Q}$ and $H^i_c(\C;\mathbb{Q}) = 0$ for $i\neq 2$, the only non-zero terms in the spectral sequence \eqref{eq1.3} when $X=\C$ are: $$	E_1^{0,2n} \cong 	E_1^{1,2n-2} \cong\mathbb{Q}.$$  We thus obtain\begin{equation}
		H^{i}_c(\uconf_n(\C);\mathbb{Q}) =\begin{cases}
			\mathbb{Q} & i=2n,2n-1\\
			0 & \text{ otherwise.}
		\end{cases}
	\end{equation} Our result, via Poincaré duality, agrees with prior computations of $H^{*}(\uconf_n(\C);\mathbb{Q})$ (see e.g. \cite{Arnold1969}, \cite{Church2012}). Of course one could have also replaced the question of computing  $H^{*}_c(\uconf_n(\C);\mathbb{Q})$ by  $H^{*}_{c,\et}(\uconf_n(\A^1);\mathbb{Q}_{\ell})$ over a field $K$, in which case the second half of Theorem \ref{thm1} gives us the desired answer.  
	
	\subsection{\small (Tuples of) polynomials with specified multiplicity of common roots}\label{sub3.2}
	In the paper \cite{Farb2015}, Farb and Wolfson studied the moduli space $\mathit{Poly}^{n,r+1}_{v}$ which, over a field $K$, they defined as \begin{align*}
		\mathit{Poly}^{n,r+1}_v:=\{(g_0, \ldots, g_{r}): g_i\in K[z] \text{ monic of degree } n, \text{ such that } {g_0,\ldots, g_r} 
		\\ \text{have no common root over } \overline{K} \text{ with multiplicity } \geq v \}
	\end{align*} When $v=1$, this is the well-known moduli space of morphisms $\P^1\to \P^r$ of degree $n$ that take $\infty\in \P^1$ to $[1:\ldots:1] \in \P^r$. Note that when $v>n$ the condition of having $v$ common roots is empty. 
	
	In \cite{Farb2015} they compute $H^*(\mathit{Poly}^{n,r+1}(\C);\Qb)$ (as well as work out some interesting arithmetic and geometric refinements via comparison theorems) by making studying the spaces over $\C$ and making use of a beautiful technique of Segal's-- that of `bringing zeroes in from infinity' (see \cite{Segal}).
	In this paper, instead of using Segal's`bringing zeroes from infinity' technique, we give a strictly algebraic proof. Let us elaborate. 
	
	We fix an algebraically closed field $K$. 
	Let $\mathit{Poly}^{n,r+1}$ be the space of all $(r+1)$-tuples polynomials of degree $n$; therefore $\mathit{Poly}^{n,r+1} \cong \A^{(n)(r+1)}$. 
	Let us define, for all $p\geq 0$, spaces (complex manifolds when $K=\C$, smooth schemes of finite type over $K$ that care actually defined over $\mathbb{Z}$): 
	\begin{align*}
		T_p=\{(z_0,\ldots, z_p), (g_0,\ldots, g_r): \text{if } z_j \text{ occurs } \lambda_j \text{ times, then  } (z-z_j)^{\lambda_jv} \text{ divide } g_i \\ \text{ for all } 0\leq i \leq r, 0\leq j\leq p\}.
	\end{align*} It immediately follows that by $T_{\bullet}$ is a symmetric semisimplicial space. In fact, by Definition \ref{def2.5} $\mathit{Poly}^{n,r+1}$ admits a symmetric semisimplicial filtration by $\A^1$ with filter gap $e=v$. Indeed, for all $p$ we have an isomorphism \begin{gather*}
		T_p \to (\A^1)^{p+1}\times \mathit{Poly}^{n-pv, r+1}\\ (z_0,\ldots, z_p), (g_0,\ldots, g_r) \mapsto (z_0,\ldots, z_p),\Bigg({g_0\over {\big((z-z_0)\ldots(z-z_p)}\big)^v}, \ldots,  {g_r\over {\big((z-z_0)\ldots(z-z_p)}\big)^v}\Bigg), 
	\end{gather*} and the face maps are finite morphisms given by \begin{gather*}
		f_i: T_p \to T_{p-1} \\ (z_0,\ldots, z_p), (g_0, \ldots, g_r) \mapsto (z_0,\ldots, \widehat{z_i}, \ldots, z_p), (g_0, \ldots, g_r) .
	\end{gather*} or equivalently, \begin{gather*}
		f_i:  (\A^1)^{p+1}\times \mathit{Poly}^{n-(p+1)v, r+1} \to  (\A^1)^{p}\times \mathit{Poly}^{n-pv, r+1}\\ (z_0, \ldots, z_p), (h_0, \ldots, h_r) \mapsto(z_0,\ldots, \widehat{z_i}, \ldots, z_p), ((z-z_i)^v h_0, \ldots, (z-z_i)^vh_r).
	\end{gather*}
	Plugging $M=\A^1$ and $X_n= \mathit{Poly}^{n,r+1}$, and in Theorem \ref{thm1} we obtain a spectral sequence which reads as:
	\begin{align}\label{eqpolynm}
		E_1^{p,q} = \big(H^q_c(\A^{p}\times \mathit{Poly}^{n-pv, r+1})\otimes \sgn_{p}\big)^{{S}_p} \implies H^{p+q}_c(\mathit{Poly}^{n,r+1}_{v}),
	\end{align}
	(where, for any $\mathbb{Z}$-scheme $S$, we mean $H^q_c(S)$ to stand for both $H^q_c(S(\mathbb{C});\mathbb{Q})$ as well as $H_{\et,c}^q(S_{/K}; \mathbb{Q}_{\ell})$,  $\ell$ coprime to $\mathrm{char}\,\, K$). Just like in the proof of Theorem \ref{thm1}, by \cite{Macdonald1962} we know that the only values of $p$ for which $E_1^{p,q}$ is nonzero $p=0, 1$; indeed \begin{equation*}
		H^*\big(((\A^1)^p)\otimes \sgn_{p}\big)^{{S}_p} \cong \begin{cases}
			H^0(\A^1) &p=1\\ 0 & \text{ otherwise.}
		\end{cases}
	\end{equation*}
	Therefore the terms in the spectral sequence \ref{eqpolynm} becomes \begin{gather*}
		E_1^{p,q} =\begin{cases}
			H_c^{2n(r+1)}(\mathit{Poly}^{n,r+1}) & (p,q) = (0, 2n(r+1))\\ H_c^{2(n-v)(r+1)+1}(\A^1\times \mathit{Poly}^{n-v,r+1}) & (p,q)= (1, 2(n-v)(r+1)+1)\\ 0 & \text{ otherwise}
		\end{cases}\\
		\cong \begin{cases}
			\Q & (p,q) = (0, 2n(r+1)), (1, 2(n-v)(r+1)+1) \\0 & \text{ otherwise.}
		\end{cases}
	\end{gather*} 
	Clearly the spectral sequence degenerates on the $E_1$ page and we obtain the following:
	\begin{cor}\label{cor7}
		Over $\C$ we have:
		\begin{equation*}
			H^i(\mathit{Poly}^{n,r+1}_{v};\Qb) \cong \begin{cases}
				\Qb & i=0, 2v(r+1)-3 \\0 &\text{ otherwise}
			\end{cases},
		\end{equation*} and we have an isomorphism of $Gal(\overline{\mathbb{F}_q}/\mathbb{F}_q)$- representations:
		\begin{equation*}
			H^i(\mathit{Poly}^{n,r+1}_{v}(\overline{\mathbb{F}_q});\Q_{\ell}) \cong \begin{cases}
				\Q_{\ell}(0) & i=0\\ \Q_{\ell}((v-n)(r+1)-1) & i= 2v(r+1)-3 \\0 &\text{ otherwise}
			\end{cases},
		\end{equation*} 
	\end{cor}
	thus recovering the cohomology part of Farb-Wolfson's \cite[Theorem 1.2]{Farb2015}, and in the special case of $v=1$ this an algebro-geometric and arithmetic analogue of Segal's \cite[Propositions 1.1 and 1.2]{Segal}.
	
	\subsection{\small Moduli space of degree $n$ morphisms $\P^1 \to \P^r$.} We continue working on the algebraically closed field $K$ fixed in the previous example. As mentioned, a special case of the previous example is that of the moduli space of degree $n$ \emph{based} maps $\P^1\to \P^r$. Now we consider the moduli space \emph{non-based} $\P^1\to \P^r$ of degree $n$ maps and prove Corollary \ref{cor4}. Note that even though Corollary \ref{cor4} follows from Theorem \ref{thm3}, which considers maps from a genus $g$ smooth projective curve for $g\geq 0$ (which will be proved in the next section), past literature supports that it's worth to work out the case $g=0$ for itself. 
	
	\begin{proof}[Proof of Corollary \ref{cor4}.]
		For $r\geq1$ define $$\Gamma_n(r) := \Big\{(s_0, \ldots, s_r): s_i\in \Gamma(\P^1,\mathscr{O}_{\P^1}(n))\Big\},$$ where $\mathcal{O}_{\P^1}(1)$ is the `hyperplane bundle' or invertible sheaf given by the sections of the universal bundle on $\P^1$, and let $$\mathcal{O}_{\P^1}(n): = \mathcal{O}_{\P^1}(1)^{\otimes n}.$$ Elements of $\Gamma(\P^1,\mathscr{O}_{\P^1}(n))$ can be thought of as homogenous polynomials of degree $n$ in two variables $x,y$; in particular, $\Gamma_n(r) \cong \A^{(r+1)(n+1)}$. 
		An element $(s_0,\ldots, s_r)\in \Gamma_n(r)$ having no common roots on $\P^1$ (which will often phrase as: $(s_0,\ldots, s_r)$ is \emph{basepoint free}, not to be confused with (non)based maps discussed above) defines a (unique, up to multiplication by $\mathbb{G}_m= K^{\times}$) map:\begin{gather*}
			\P^1\to\P^r\nonumber \\
			z\mapsto [s_0(z):\ldots: s_r(z)].
		\end{gather*} Conversely, any map $\P^1\to \P^r$ of degree $n$ is given by a basepoint free $(r+1)$-tuple of sections $(s_0, \ldots, s_r)\in \Gamma_n(r)$. We say that an element $[s_0: \ldots : s_r]\in \P\Gamma_{n}(r)$ is basepoint free if $(s_0,\ldots, s_r)\in \Gamma_n(r)$ is. Let the locus of the basepoint free elements of $\P\Gamma_n(r)$ be denoted by $\mathrm{Mor}_n(\P^1,\P^r);$ i.e. $$\mathrm{Mor}_n(\P^1,\P^r);:= \{[s_0: \ldots: s_r]: \nexists [a:b]\in \P^1 \text{ such that } s_i([a:b])=0  \text{ for all }i\}; $$ 
		parts of literature also call it the \emph{Hurwitz space of degree $n$ morphisms $\P^1\to \P^r$.} Let $Z_n(r): = \Gamma_n(r)- \mathrm{Mor}_n(\P^1,\P^r);$  be the \emph{discriminant locus} i.e. $$Z_n(r):= \{[s_0:\ldots: s_r]: \exists [a:b]\in \P^1 \text{ such that } s_i([a:b])=0  \text{ for all }i\}.$$ Note that $Z_n(r)$ is defined over $\mathbb{Z}$; indeed, it's cut out by polynomials in $\P\Gamma_n(r)$ defined over $\mathbb{Z}$. In turn $U_n(r)$ is defined over $\mathbb{Z}$.
		
		We apply Theorem \ref{thm2} to compute $H^*(\mathrm{Mor}_n(\P^1,\P^r))$ by plugging in $X_n:=\P \Gamma_n(r)$ and showing it admits a semisimplicial filtration by $\P^1$, with the space of topological $\P^1$- indecomposables being $\mathrm{Mor}_n(\P^1,\P^r)$ and the filter gap $e=1$. 
		To begin, note that for each $p$ we have face maps given by \emph{adding a basepoint}: \begin{gather}
			f_{i} : (\P^1)^{p+1} \times X_{n-(p+1)} \to (\P^1)^{p} \times X_{n-p} \nonumber\\
			\big([a_0:b_0],\ldots, [a_p:b_p]\big),\,\, [s_0:\ldots: s_r] \mapsto \nonumber \\ \big([a_0,b_0],\ldots, \widehat{[a_i:b_i]},\ldots [a_p:b_p]\big), \,\, [(b_ix-a_iy)s_0:\ldots: (b_ix-a_iy)s_r]\label{eq3.5}
		\end{gather} In other words, the hypercover under consideration is the following:
		\begin{gather*}
			\cdots\cdots (\P^1)^3\times \P\Gamma_{n-3}(r)\mathrel{\substack{\textstyle\rightarrow\\[-0.5ex]
					\textstyle\rightarrow \\[-0.5ex]
					\textstyle\rightarrow}} 	(\P^1)^2 \times \P\Gamma_{n-2}(r) \rightrightarrows \P^1\times \P\Gamma_{n-1}(r)\to \P\Gamma_{n}(r)
		\end{gather*} with the unlabelled arrows denoting the face maps $f_i$.
		It is almost immediate that the face maps satisfy all the conditions from Definition \ref{def2.10} with $M=\P^1$ and $X_n= \P\Gamma_n(r)$. Plugging them in Theorem \ref{thm2}, we obtain a second quadrant spectral sequence which reads as \begin{gather*}
			E_1^{-p,q} = \begin{cases}
				H^q(\P\Gamma_{n}(r))(0) & p=0,\\ H^{q-2r}(\P^1\times \P\Gamma_{n-1}(r))(-1) & p=1,\\ H^0(\P^1)\otimes H^2(\P^1)\otimes H^{q-4r-2}(\P\Gamma_{n-2}(r))(-2) & p=2,\\ 0 & \text{ otherwise},
			\end{cases}
		\end{gather*} with the differentials given by the alternating sum of the Gysin pushforwards induced by the face maps, which is what we shall compute now. 
		\begin{itemize}
			\item \textit{Computing $d_1^{1,q}: E_1^{-1,q} \to E_1^{0,q}$.} 
			
			For simplicity we denote the differential by $d_1^1$. Let $$\iota: \P\Gamma_{n-1}(r) \hookrightarrow \P\Gamma_{n}(r)$$ denote the inclusion given by adding a basepoint. Choose generators $\1\in H^0(\P^1)$ and $e\in H^2(\P^1)$, and let $h$ denote the hyperplane class in $\P\Gamma_n(r)$. Then we claim that: \begin{gather*}
				d_1^1={f_0}_*: H^{*-2r}(\P^1\times \P\Gamma_{n-1}(r)) \to   H^*(\P\Gamma_{n}(r))\\ \1\otimes \iota^*\alpha + e\otimes \iota^*\alpha' \mapsto \alpha h^r+ \alpha' h^{r+1}
			\end{gather*} is a map of $H^*(\P \Gamma_{n}(r))$-modules, where $\alpha, \alpha' \in H^*(\P\Gamma_{n}(r))$. To see this, first note that $$\iota^*: H^*(\P\Gamma_{n}(r)) \to H^*(\P\Gamma_{n-1}(r))$$ is a surjection; next, the image of the fundamental class $$[\P^1\times \P\Gamma_{n-1}(r)]\in H^0(\P^1\times \P\Gamma_{n-1}(r))$$ is the locus of elements in $\P\Gamma_n(r)$ that has a basepoint i.e. $Z_n(r)$, which is rationally equivalent, and thus cohomologous, to (a multiple of) $h^r$; and finally, for a fixed point $[a:b]\in \P^1$, the locus given by $$\{[s_0:\ldots: s_r]\in \P\Gamma_{n}(r): s_i([a:b])=0\}$$ is rationally equivalent, and in turn cohomologous, to (a multiple of) $h^{r+1}$. For the sake of simplicity we won't bother ourselves with the scalar multiples, which is fine because we're working over $\Q$.
			
			The Gysin pushforward $d_1^1={f_0}_*$ surjects onto the ideal generated by $h^r$ in $H^*(\P \Gamma_{n}(r))$. Indeed, the preimage of $h^{r+i}$ is given by \begin{gather*}
				d_1^1(\1\otimes \iota^* h^i) =h^{r+i} = d_1^1(e\otimes \iota^*h^{i-1}) \text{ for } i\geq 1\\ d_1^1([\P^1\times \P\Gamma_{n-1}(r)]) =h^r,
			\end{gather*} which shows that the image of $d_1^1$ is the ideal generated by $h^r$ in $H^*(\P\Gamma_{n}(r))$. The kernel of $d_1^1$ is given by elements of the form $(h-e)\otimes \iota^*(\alpha)$ for all $\alpha \in H^*(\P\Gamma_{n}(r))$. Again, recalling that $\iota^*: H^*(\P\Gamma_{n}(r)) \to H^*(\P\Gamma_{n-1}(r))$ is a surjection, we conclude that $\textrm{Kernel}(d^1_1) $ is generated by elements of the form $$(h-e)\otimes \beta,\,\,\,\,\,\, \beta\in H^*(\P\Gamma_{n-1}(r)).$$
			
			The upshot is that on the $E_2$ page, for $p=0$ we have:\begin{gather}\label{eq3.7}
				E_2^{0,q} = \begin{cases}
					\Q(0) & q=2j, \,\, 0\leq j\leq 2(r-1)\\0 & \text{ otherwise.}
				\end{cases}
			\end{gather} 
			
			\item  \textit{Computing $d_1^{2,q}: E_1^{-2,q} \to E_1^{-1,q}$.} 
			
			For simplicity, we denote the differential by $d_1^2$. Like before, let $$\iota: \P\Gamma_{n-2}(r) \hookrightarrow \P\Gamma_{n-1}(r)$$ denote the inclusion given by adding a basepoint, and  let $h$ denote the hyperplane class in $\P\Gamma_{n-1}(r)$. Let us also keep in mind, like before, that  $$\iota^*: H^*(\P\Gamma_{n-1}(r)) \to H^*(\P\Gamma_{n-2}(r))$$ is a surjection. Then the way we computed ${f_0}_*$ above works verbatim, and we have \begin{gather*}
				{f_0}_*:H^0(\P^1)\otimes H^2(\P^1)\otimes H^{*-2r-2}(\P\Gamma_{n-2}(r)) \to H^{*}(\P^1\times \Gamma_{n-1}(r))\\
				\1 \otimes e\otimes \iota^*\alpha \mapsto e\otimes\alpha h^r 
			\end{gather*} and 
			\begin{gather*}
				{f_1}_*:H^0(\P^1)\otimes H^2(\P^1)\otimes H^{*-2r-2}(\P\Gamma_{n-2}(r)) \to H^{*}(\P^1\times \P\Gamma_{n-1}(r))\\
				\1 \otimes e\otimes\iota^* \alpha \mapsto \1\otimes\alpha h^{r+1},
			\end{gather*} and therefore$$d_1^2(1\otimes e \otimes \iota^*\alpha)= 1\otimes \alpha h^{r+1}- e\otimes \alpha h^r = (h-e)\otimes \alpha h^r.$$ Note that $d_1^2$ is injective, and the image is generated by $h^r$ in $H^*(\P^1\times \P\Gamma_{n-1}(r))$. Consequently, on the $E_2$ page we have: \begin{flalign*}
				&E_2^{-1,q} = \begin{cases}
					\Q(-r) & p=1, q=2j+2r+2, \,\, 0\leq j\leq 2(r-1)\\0 &\text{ otherwise}
				\end{cases},
				\\ & E_2^{-2,q} = 0,  \text{ for all } q.
			\end{flalign*} 
		\end{itemize}
		In effect on the $E_2$ page all differentials vanish; the spectral sequence degenerates and we obtain $$H^*(\mathrm{Mor}_n(\P^1,\P^r);\Q)\cong {\Q[h]\over h^r} \otimes \wedge\Q\{t\}$$ where $h$ has cohomological degree $2$, and $t$ (which corresponds to $e-h\in \textrm{Ker}(d_1^1)$) has cohomological degree $2r+1$. Furthermore, over a field $\kappa$, with algebraic closure $\overline{\kappa}$, we have an isomorphism of $Gal (\overline{\kappa}/\kappa)$-representations:
		
		\begin{gather*}
			H^i_{\et}(U_n(r);\Qb_{\ell}) = \begin{cases}
				\Qb_{\ell}(-j) & i=2j, 0\leq j\leq r-1\\ \Qb_{\ell}(-(j+1)) & i=2j+1, r\leq j\leq 2r-1\\ 0 & \text{ otherwise.} 
			\end{cases}
		\end{gather*}
		This completes our proof of Corollary \ref{cor4}.
	\end{proof}
	
	\subsection{\small Moduli space of degree $n$ morphisms $C \to \P^r$, $g(C)\geq 0$.} When $g=0$, we have $C \cong \P^1$ and we discussed it above. Now, let $C$ be a fixed smooth projective curve of genus $g$ where $g\geq 0$, and fix a positive integer $r$. We compute the (stable) cohomology of the moduli space $\mathrm{Mor}_n(C,\P^r)$ of degree $n$ morphisms $C\to \P^r$. 
	
	A degree $n$ morphism $C\to\P^r$ is equivalent to the following data:\begin{itemize}
		\item a line bundle $L$ of degree $n$ on $C$,
		\item an $(r+1)$-tuple $(s_0,\ldots, s_r)$ where $s_i\in H^0(C,L)$ 
		\item the sections $s_0,\ldots, s_r$ satisfy the condition that they have no common zeroes (also known as $\{s_0,\ldots, s_r\}$ is \emph{basepoint free}).
	\end{itemize} Then, $\mathrm{Mor}_n(C,\P^r)$ is a Zariski open dense subset of the smooth projective variety $X_n$ defined by \begin{gather}\label{eq3.8}
		X_n:=\{L,[s_0:\ldots:s_r]: L \in \pic^n(C), s_i\in H^0(C,L) \text{ for all }i\}.
	\end{gather}When $n\geq 2g$ (for $g\geq 2$, even $n\geq 2g-1$ works for our purposes), by the Riemann-Roch theorem $\dim H^0(C,L) =n-g+1$ for all $L\in \pic^n(C)$, and $X_n$, in turn, is isomorphic to the projectivization of a vector bundle $E_n$ on $\pic^n(C)$ whose fibres are isomorphic to $\A^{(n-g+1)(r+1)}$. To elaborate, let $P(n)$ be a \emph{Poincaré line bundle for $C$ of degree $n$} (see \cite[Chapter IV, Section 2]{ACGH} for the definition of a Poincaré line bundle and its properties) and let $$\nu_n:C\times \pic^n(C)\to \pic^n(C)$$ be the projection to the second factor. Then for each $n\geq 2g-1$ we have a vector bundle $$E_n={\nu_n}_*P(n) \to \pic^n(C),$$ with the fibre over a point $[L]\in  \pic^n(C)$ being $$\big(E_n\big)_L =H^0(C,L)\cong \A^{n-g+1},$$ and then $X_n$ can be equivalently described as: $$\rho_n: X_n = \P E_n^{r+1}\to \pic^n(C)$$ where $E_n^{r+1}$ is the $(r+1)$-fold fibre product of $E_n$ over $\pic^n(C)$. We show that $X_n$ admits symmetric semisimplicial filtration by $C$ with filter gap $e=1$, and use Theorem \ref{thm2} to compute $H^i(\morn;\Qb)$ for $i\leq n-2g+1$.
	
	Some notations before we start proving Theorem \ref{thm3}: we suppress the coefficient field and just write $H^*(X)$ to stand for $H^*(X;\Qb)$ until we come to the point where we have keep track of weights, and in particular, the necessary Tate twists.
	
	\textsc{Proof outline:} Because the proof is very involved and computational, we split the proof into several parts which we outline before we begin the proof.
	\begin{enumerate}[(i)]
		\item We show $X_n$, as defined in \eqref{eq3.8}, admits symmetric semisimplicial filtration by $C$ with filter gap $1$. 
		\item We plug in $M=C$ and $X_n$ in the spectral sequence \eqref{eq3.8} from Theorem \ref{thm2}, and compute the terms $E_1^{-p,q}$ of the (second-quadrant) spectral sequence for $0\leq p\leq n-2g$.
		\item We compute the differentials on the $E_1$ page in the range $0\leq p\leq n-2g+1$ and deduce the $E_2$ terms.
		\item We show that $E_2^{-p,q} =  E_{\infty}^{-p,q}$ for $0\leq p\leq n-2g$.
	\end{enumerate}
	
	\begin{proof}[Proof of Theorem \ref{thm3}]\begin{enumerate}[Step 1.]
			
			\item 	\textit{Show $X_n$ admits symmetric semisimplicial filtration.}
			
			To this end, observe that an equivalent description of $\P E_n$ is that it is the space of all effective divisors on $C$ of degree $n$. Indeed, the fibre of the map $\P E_n\to \pic^n(C)$ over $\mathcal{O}_C(D)\in \pic^n(C)$ is the complete linear system of all effective divisors $D'$ of degree $n$ that are rationally equivalent to $D$ (often written as $D'\sim D$), and $$\{D': D' \text{  is effective of degree }n, D'\sim D\} = \P H^0(C,\mathcal{O}(D)).$$ In turn, for each $x\in C$, we have a commutative diagram: 
			
			\[
			\begin{tikzcd}
				\P E_{n} \arrow{r}{t_x^{\mathit{eff}}} \arrow[swap]{d} & \P E_{n+1} \arrow{d} \\
				\mathit{Pic}_n(C) \arrow{r}{t_x} &  \mathit{Pic}_{n+1}(C)
			\end{tikzcd}
			\]
			where \begin{gather*}
				\txf: C\times \P E_n\to \P E_{n+1}\\ x, D\mapsto x+D
			\end{gather*} is the map of \emph{adding a point $x$} on effective divisors, and 
			\begin{gather*}
				t_x:\pic^n(C)\xrightarrow{\cong} \pic^{n+1}(C) \\ x, \mathcal{O}_C(D) \mapsto \mathcal{O}_C(x+D)
			\end{gather*} is the \emph{translation by $x$} map on the Picard group, which is naturally an isomorphism. Now observe that $\txf$ is a relative linear embedding of $\P E_n$ in $\P E_{n+1}$ as schemes over $\pic^{n+1}(C)\cong t_x(\pic^n(C))$. This is because of the following. A Poincaré bundle $P(n)$ is $\nu_n$-relatively very ample when $n\geq 2g$ because it is fibrewise very ample for the proper map $$\nu_n: C\times \pic^n(C) \to \pic^{n+1}(C)$$ (see \cite[Chapter 1, Section 1.7]{Lazarsfeld2004} or \cite[Section 4.7.1]{Grothendieck1961}). Therefore the relative evaluation map of locally free sheaves on $C\times\pic^n (C)$: $$\mathit{ev}_{{x}\times \pic^n(C)}: \nu_n^*{\nu_n}_*P(n) \to \mathcal{O}_{C\times \pic^n C}$$ is surjective and the kernel, which is a locally free sheaf, is a relative hyperplane bundle in $\P E_{n+1}$ and is the image $\txf(\P E_n)$ by definition. All this is to conclude that the addition by $x$ map on the space of effective divisors has a natural `lift' to a map of adding a point $x$ on the vector bundle $E_n$: $$t_x^{\mathit{glob}}:E_n\to E_{n+1}$$ and in turn it results in a similar addition by $x$ map on $\P E^{r+1}_n$ as follows: \[
			\begin{tikzcd}
				\P E^{r+1}_{n} \arrow{r}{t_x^{r}} \arrow[swap]{d} & \P E^{r+1}_{n+1} \arrow{d} \\
				\mathit{Pic}_n(C) \arrow{r}{t_x} &  \mathit{Pic}_{n+1}(C)
			\end{tikzcd}
			\] where we have \begin{gather*}
				t^r_x: \P E^{r+1}_{n} \to \P E^{r+1}_{n+1}\\ [s_0:\ldots: s_r] \mapsto [t_x^{\mathit{glob}}(s_0):\ldots: t_x^{\mathit{glob}}(s_r)].
			\end{gather*}
			Going through the whole drill above for all $x \in C$, one gets a natural addition map  \[
			\begin{tikzcd}
				C\times \P E^{r+1}_{n} \arrow{r}{A} \arrow[swap]{d} & \P E^{r+1}_{n+1} \arrow{d} \\
				C\times	\mathit{Pic}_n(C) \arrow{r}{A^{\mathit{rat}}} &  \mathit{Pic}_{n+1}(C)
			\end{tikzcd}
			\]
			where \begin{gather*}
				A:C\times \P E^{r+1}_{n} \to \P E^{r+1}_{n+1} \\ x, \big(L, [s_0:\ldots: s_r]\big) \mapsto \Big(L\otimes \mathcal{O}_C(x), t_x^r\big([s_0:\ldots: s_r]\big)\Big)
			\end{gather*} is, just like in the case of $C=\P^1$ in the previous example, \emph{adding a basepoint}, and where $$A^{\mathit{rat}}: C\times \pic^n(C)\to \pic^{n+1}(C)$$ is the addition map on rational equivalence classes of divisors i.e. $$A^{\mathit{rat}}(x,\mathcal{O}_C(D)) = \mathcal{O}_C(x+D)$$ and the image $A(C\times \P E^r_n)$ is precisely given by $$\big\{L,[s_0:\ldots, s_r]: L\in \pic^{n+1}(C), s_i\in H^0(C,L), s_0,\ldots, s_r \text{ have a common zero}\big\}.$$ Note that the adding a basepoint maps, even though defined set-theoretically, are radiciel maps (they are injective on the $\C$-points) and in fact it is easy to check that they are closed embeddings.
			
			The rest of the proof essentially follows that of the case of $C=\P^1$ from the previous example. Define a semisimplicial space by $$T_{p-1}:= C^p\times \P E^{r+1}_{n-p}$$ for $p\geq 1$, where $E_{n-p}$ is a vector bundle as long as $n-p\geq 2g$, and where the face maps are given by adding a basepoint: \begin{gather*}
				f_i: C^{p+1}\times \P E^{r+1}_{n-(p+1)} \to C^{p}\times \P E^{r+1}_{n-p}\\
				(x_0,\ldots, x_p),\big(L, [s_0:\ldots: s_r]\big) \mapsto  (x_0,\ldots, \hat{x_i}, \ldots, x_p),\big(L\otimes \mathcal{O}_C(x_i), t_x^r[s_0:\ldots: s_r]\big) 
			\end{gather*} i.e. the $i^{th}$ face map is just the map $A$ using the $i^{th}$ copy of $C$, with identity on the remaining copies of $C$, and in fact $f_0=A$. The semisimplicial space thus defined is clearly symmetric semisimplicial. For $p=0$, we define $T_{-1}:= X_n$.
			
			\item \textit{Apply Theorem \ref{thm2} and compute the $E_1$ terms.}
			Plugging $M=C$, $X_n = \P E^{r+1}_{n}$ and $e=1$ in Theorem \ref{thm2}, we obtain a second quadrant spectral sequence which reads as:\begin{gather}\label{eq3.9}
				E_1^{-p,q}= H^{q-2pr}\big((C\times \P E^{r+1}_{n-p})\otimes \sgn_{p}\big)^{{S}_p}(-pr) \implies H^{q+p}(\morn).
			\end{gather}
			where we keep a record of the Tate twists (given by, for each $p$, the codimension of $T_{p-1}$ in the geometric realization of $T_{\bullet}$) to keep track of the Hodge structures, and where the differentials are given by the alternating sum of the Gysin pushforwards induced by the face maps. To this end note that by \cite{Macdonald1962} we know that: \begin{flalign}\label{eq3.10}
				\big(H^*(C)^{\otimes p}\otimes \sgn_{p}\big)^{{S}_p}  & \cong H^0(C)\otimes \Sym^{p-1}H^1(C) \nonumber \\ & \bigoplus H^2(C)\otimes \Sym^{p-1}H^1(C) \nonumber \\
				& \bigoplus H^0(C)\otimes H^2(C)\otimes \Sym^{p-2}H^1(C) \nonumber\\ &\bigoplus \Sym^p H^1(C).
			\end{flalign}
			Therefore to have a complete understanding of $H^*(\P E^{r+1}_n)$ for all $n\geq 2g$ we need to know the Chern classes of $E^{r+1}_{n}$. For $r=0$ we have $E_n^1=E_n$ and the Chern classes of $E_n$ can be computed for example, directly using Grothendieck-Riemann-Roch, or via ad-hoc methods to give us $$c_i(E_n) = (-1)^i {\theta^i \over i!}\,\,\, i=0, \ldots, g$$
			where $\theta$ is the fundamental class of the theta divisor (several proofs are available in \cite[Sections 4, 5, Chapter VII and Section 1, Chapter VIII]{ACGH}). Using the Whitney sum formula we obtain the Chern classes of $E^{r+1}_n$: \begin{align*}
				c_i(E^{r+1}_n)= \sum_{\substack{0\leq i_0,\ldots, i_r\leq g\\0.i_0+ 1.i_1+2.i_2+\ldots + r i_r =i}} (-1)^{i} {\theta^i \over  i_0!\ldots i_r!} \\ = (-1)^i {\binom{r+i}{i}}{\theta^i\over i!}.
			\end{align*} In turn, let $N_0:= (n-g+1)(r+1)$, the dimension of the fibres of $$E^{r+1}_n\to \pic^n(C),$$  and let $h$ denote the relative hyperplane class i.e. $$h= c_1(\mathcal{O}_{\rho_n}(1)) \in H^2(\P E^{r+1}_n),$$ then $H^*(\P E^{r+1}_n)$, which is an algebra on $$H^*(\pic^n (C))\cong\wedge(H^1(C)),$$ is given by \begin{gather}\label{eq3.11}
				H^*(\P E^{r+1}_n)\cong {H^*(\pic^n (C))[h]\over {h^{N_0}+\rho_n^*c_1(E^{r+1}_n) h^{N_0-1}+ \ldots + \rho_n^*c_g(E^{r+1}_n)h^{N_0-g}}}.
			\end{gather}	
			Let $p$ be such that $n-p\geq 2g$ and let $$N_p:= (n-p-g+1)(r+1) = N_0-p(r+1),$$ the dimension of the fibres of $E^{r+1}_{n-p}\to \pic^{n-p}(C)$, then combing \eqref{eq3.10} and \eqref{eq3.11} we have a complete description of the $E_1$ terms of the spectral sequence above. We note here that since $n-p\geq 2g$,  we have that $$N_p-g = (n-p-g+1)(r+1)-g> r.$$ This observation will be useful later.
			
			\item   \textit{Computing the differentials $d_1^p: E_1^{-p,*}\to E_1^{-(p-1),*+2r}$.}
			
			Following previously introduced notations, let $h= c_1(\mathcal{O}_{\rho_{n}}(1))$, and for all $p$ satisfying $n-p\geq 2g$, let $$\iota: \P E^{r+1}_{n-p}\to \P E_{n-(p+1)}^{r+1}$$ denote the closed embedding induced by adding a basepoint $x$ (an abuse of notation that won't cause any confusion down the way) i.e. $\iota = t_x^r$ for a chosen $x\in C$. Note that $\iota$ is fibrewise a linear embedding, up to translation of $\pic^{n-p}(C)$ by $x$. Finally, let $e\in H^2(C)$ be the class of a point, $\1$ the fundamental class of $C$, and let $c_1,\ldots, c_{2g}$ be the standard basis of $H^1(C)$ and because $H^*(\pic^n(C)) \cong \wedge H^1(C)$, let $\overline{c_1},\ldots, \overline{c_{2g}}$ be the image of $c_1, \ldots, c_{2g}$ under the aforementioned isomorphism. 
			
			First, we observe that \begin{align*}
				d_1^{1}: H^*(C\times \P E^{r+1}_{n-1}) \to H^*(\P E^{r+1}_n)\\ [C\times \P E^{r+1}_{n-1}] \mapsto h^r\\ e\mapsto h^{r+1}\\ c_i\mapsto \overline{c_i}h^r, & \text{ for all }i.
			\end{align*} is a map of $H^*(\P E_{n}^{r+1})$-modules, and in turn $$\iota^*\alpha + e \iota^*\beta + \sum_{i=1}^{2g} c_i \iota^* \gamma_i \xmapsto{d_1^1} \alpha h^r+ \beta h^{r+1}+ \sum_{i=1}^{2g}  \overline{c_i}\gamma_i h^r,$$ where $\alpha, \beta, \gamma_1,\ldots, \gamma_{2g} \in  H^*(\P E^{r+1}_n)$.
			Indeed, the justification for the formula for $d_1^1$ in the previous case of $C=\P^1$ holds almost verbatim here. We know $$\iota^*: H^*(\P E^{r+1}_{n}) \to H^*(\P E^{r+1}_{n-1})$$ is a surjection; next, for a fixed point $x\in C$, the image $t_x^r(\P E^{r+1}_{n-1}$) is rationally equivalent, and in turn cohomologous, to (a multiple of) $h^{r+1}$, and finally, that the image of the fundamental class $[C\times \P E_{n-1}^{r+1}]\in H^0(C \times \P E_{n-1}^{r+1})$ is rationally equivalent, and thus cohomologous, to (a scalar multiple of) $h^r$, can be seen as in the following way. Recall that a Poincaré bundle $P(n)$ is $\nu_n$-relatively very ample for all $n\geq 2g$, which in turn induces a relative embedding of $C\times \pic^n(C)\xrightarrow{i_n}\P E_n$ over $\pic^n(C)$ i.e. 
			\[
			\begin{tikzcd}
				C\times \pic^n(C) \arrow{rr}{i_n} \arrow[swap]{dr}{\nu_n} & & \P E_n \arrow{dl}{} \\[10pt]
				& \mathit{Pic}^n(C)
			\end{tikzcd}
			\] which, over a point $[L]\in \pic^n(C)$ is merely an embedding of $C\hookrightarrow \P(H^0(C,L)^*)$ under the complete linear system of $L$. Now, $ \P E_n $ is linearly embedded in $\P E^{r+1}_n$ over $\pic^n(C)$, and let $i_n$ still denote the composition $$C\times \pic^n(C) \xhookrightarrow{i_n}  \P E_n  \hookrightarrow \P E^{r+1}_n, $$which makes $i_n(C\times \pic^n(C))$ in $\P E^{r+1}_n$ homologous to (a scalar multiple of) the Poincaré dual of $h \in H^2(\P E^{r+1}_n)$  In turn, the image of the $[C\times \P E_{n-1}^{r+1}]$ under the Gysin map ${f_0}_*$ is given by \begin{align*}
				{f_0}_*\big([C\times \P E_{n-1}^{r+1}]\big) =h^{r+1}\frown i_n(C\times \pic^n(C))\\ =h^r.
			\end{align*}
			Yet again, for the sake of simplicity we won't bother ourselves with the scalar multiples, which is fine because cohomology with $\Qb$ coefficients.
			Noting that $$\overline{c_i}(e-h)+ h(c_i-\overline{c_i}) = \overline{c_i}e - \overline{c_i}h +hc_i -h\overline{c_i} =\overline{c_i} e- c_i h,$$ it is now easy to check that the kernel of $d_1^1$ is given by: \begin{align*}
				H^*(\P E^{r+1}_{n-1})(e-h) [2r] \bigoplus _{1\leq i\leq 2g}	H^*(\P E^{r+1}_{n-1})(c_i-\overline{c_i})[2r], &&(i= 1, \ldots, 2g)
			\end{align*} where $[2r]$ denotes a shift in the cohomological degree by $2r$, and which is viewed as a $\iota^* H^*(\P E_{n}^{r+1})\cong H^*(\P E_{n-1}^{r+1})$-module. The cokernel of $d^1_1$, which forms $E_2^{0,*}$ is given by  $$H^*(\pic^n(C))[h]\over h^r$$(where note that, as observed before $r<N_0-g$, see \eqref{eq3.11}).
			
			Now we work out the differential for $p=2$ by computing the Gysin pushforwards by each of the face maps: \begin{gather*}
				{f_0}_*(\1\otimes e) = e h^r, \,\,\,\, 	{f_1}_*(\1\otimes e) = h^{r+1}  \implies d_1^2 (\1\otimes e) =  (e-h) h^{r}, \\ 	{f_0}_*(e\otimes c_i) = c_i h^{r+1}, \,\,\,\, 	{f_1}_*(e\otimes c_i) = e \overline{c_i}h^r \implies d^2_1(e\otimes c_i)= (c_ih- e\overline{c_i}) h^r,\\
				{f_0}_*(\1\otimes c_i) = c_i h^r, \,\,\,\, 	{f_1}_*(\1\otimes c_i) = \overline{c_i} h^r \implies d^2_1(\1\otimes c_i) =  (c_i -  \overline{c_i}) h^r , \\ 	d_1^2(c_ic_j)=0,
			\end{gather*} where the last equality follows form the fact that on $\Sym^p H^1(C)$ for $p\geq 2$, the alternating sum of face maps is, by definition, $0$.
			Recalling our earlier remark that $r<N_1-g$, we see that the $E_2^{-1,*}$ terms, as an $H^*(\P E_{n-2}^{r+1})$-module, are given by:
			\begin{align*}
				{H^*(\pic^{n-1}(C))(-r)[h]\over h^r}(e-h)[2r] \\ \bigoplus_{1\leq i\leq 2g} {H^*(\pic^{n-1}(C))(-r)[h]\over h^r}(c_i-\overline{c_i})[2r].
			\end{align*}Whereas the kernel of $d_1^2$ is generated by exactly what one expects: as a $H^*(\P E^{r+1}_{n-2})$-module, we have  \begin{align*}
				\mathrm{Ker}(d_1^2)=  \bigoplus_{1\leq i\leq 2g} H^*(\P E^{r+1}_{n-2})\big(e\otimes c_i  - 1\otimes c_i h +\1\otimes e \overline{c_i}\big)[4r]\\ \bigoplus_{1\leq i,j\leq 2g} H^*(\P E^{r+1}_{n-2})(c_ic_j)[4r].
			\end{align*}
			For $p= 3$ we have $d_1^{3}: E_1^{-3,*}\to E_1^{-2,*}$ given by:
			\begin{gather*}
				d_1^3 (\1\otimes e \otimes c_i) =   e\otimes c_i h^r- \1\otimes c_i h^{r+1} +\1\otimes e \overline{c_i} h^{r}\impliedby	\begin{cases}
					{f_0}_*(\1\otimes e \otimes c_i) = e\otimes c_i h^r,\\ {f_1}_*(\1\otimes e\otimes c_i) = \1\otimes c_i h^{r+1}\\ {f_2}_*(\1\otimes e\otimes c_i) = \1\otimes e \overline{c_i} h^{r}
				\end{cases}\\
				d_1^3(e\otimes c_ic_j) = c_ic_j h^{r+1}, \\
				d_1^3(\1\otimes c_ic_j) = c_ic_j h^r, \\ 	d_1^3(c_ic_jc_k)=0,
			\end{gather*} where, for the last three equalities, recall again that on $\Sym^p H^1(C)$ for $p\geq 2$, the alternating sum of face maps is, by definition, $0$. Therefore the $E_1^{-2,*}$ terms defined by $\mathrm{Ker}(d^2_1)/ \mathrm{Image}(d^3_1)$ is given by: 
			\begin{align*}
				\bigoplus_{1\leq i\leq 2g}{H^*(\pic^{n-2}(C);\Q(-2r))[h]\over h^r}\big(e\otimes c_i  - 1\otimes c_i h +\1\otimes e \overline{c_i}\big)[4r] \\ \bigoplus_{1\leq i,j\leq 2g}{H^*(\pic^{n-2}(C);\Q(-2r))[h]\over h^r}(c_ic_j)[4r].
			\end{align*} 
			The formula for the differentials in the case of $p\geq 3$ mimics that of $p=3$, and we have:
			
			\begin{gather*}
				\1\otimes e\otimes c_1\ldots c_{p-2} \mapsto \Big((e\otimes c_1\ldots c_{p-2}) - (\1 \otimes c_1\ldots c_{p-2}) h\Big)h^r,\\
				e\otimes c_1\ldots c_{p-1} \mapsto c_1\ldots c_{p-1} h^{r+1},\\  \1 \otimes c_1\ldots c_{p-1} \mapsto c_1\ldots c_{p-1} h^{r} \\ c_1\ldots c_p \mapsto 0
			\end{gather*}
			It is now easy to check that  \begin{gather*}	\mathrm{Ker}(d_1^p)/\mathrm{Image} (d_1^{p+1})
				\\ = \bigoplus_{1\leq i\leq 2g}{H^*(\pic^{n-p}(C))(-pr)[h]\over h^r}\big(e\otimes c_1\ldots c_{p-1} - \1 \otimes c_1\ldots c_{p-1} \big)[2pr]\\ \bigoplus_{1\leq i,j\leq 2g}{H^*(\pic^{n-p}(C))(-pr)[h]\over h^r}(c_1\ldots c_p)[2pr].
			\end{gather*}
			
			\item \textit{Analysing the $E_2$ page to show $E_2^{-p,q} = E_2^{-p,\infty}$ for $0\leq p\leq n-2g$.}
			
			That the differentials on the $E_2$ page vanish for $p\leq n-2g$ follow simply from weight considerations- the space $T_{\bullet}$ consists of smooth projective varieties and thus their $n^{th}$ cohomology is pure of weight $n$. Now observe the following: from Lemma \ref{lem2.11} we have an equality $$\mathrm{R}\Gamma_c(\P E^{r+1}_n,C^{\bullet}(\Qb_{\P E^{r+1}_n})) = \mathrm{R}\Gamma_c (\P E^{r+1}_n, j_{!}\Qb_{\morn})$$ in the derived category of constructible sheaves over $\P E^{r+1}_n$, where $C^{\bullet}(\Qb_{\P E^{r+1}_n})$ denotes complex in \eqref{eq2.8}, with $X_n:= \P E^{r+1}_n$; on the other hand, for any $N\in \mathbb{N}$ we have $$\mathrm{R}^i\Gamma_c(\P E^{r+1}_n, C^{\bullet}(\Qb_{\P E^{r+1}_n})) \cong \mathrm{R}^i\Gamma_c(\P E^{r+1}_n, C^{\bullet}(\Qb_{\P E^{r+1}_n})/\tau_{\geq N} C^{\bullet}(\Qb_{\P E^{r+1}_n})$$ for all $i\geq 2(r+1)-2N$, where $\tau_{\geq N} C^{\bullet}(\Qb_{\P E^{r+1}_n})$ denotes the truncated complex up to the $(N-1)$ term and this is because $\tau_{\geq N} C^{\bullet}(\Qb_{\P E^{r+1}_n})$ is supported on complex codimension $N$ in $\P E^{r+1}_n$. Therefore the cohomology of $\morn$ up to degree $n-2g$ is solely dictated by the $E_2$ page. 	To this end, let $$t:= (e-h)$$ which has degree $(-1, 2r+2)$ and let $$\alpha_i:= c_i-\overline{c_i},\,\,\,\,\, i=1,\ldots, 2g$$ which has degree $(-1,2r+1)$. Clearly for $3\leq p\leq n-2g$, the element $t\alpha_{i_1}\ldots \alpha_{i_p}$, which is of degree $(-(p+1), 2r+2+p(2r+1))$, when expanded, gives us  \begin{flalign*}
				t\alpha_{i_1}\ldots \alpha_{i_p}\\ &= (e-h)(c_{i_1}-\overline{c_{i_1}})\ldots (c_{i_p}-\overline{c_{i_p}}) \\&= (e-h)\prod_{j=1}^{p}c_{i_j} +\Big\{ \text{lower order terms as a polynomial on } c_{i_1},\ldots, c_{i_p}\Big\} \\ & = (e-h)\prod_{j=1}^{p}c_{i_j}
			\end{flalign*} because the lower order terms are all $0$ in $$\big(H^2(C)\oplus H^0(C)\big)\bigotimes \Sym^p H^1(C)\otimes H^*(\pic^{n-(p+1)}(C))[h]/h^r,$$ thanks to the alternating action of ${S}_{p+1}$. 
			
			On the other hand $\alpha_{i_1}\ldots \alpha_{i_{p+1}}$, which is of degree $(-(p+1), (p+1)(2r+1))$,  when expanded, gives us
			\begin{flalign*}
				\alpha_{i_1}\ldots \alpha_{i_{p+1}}\\ &= (c_{i_1}-\overline{c_{i_1}})\ldots (c_{i_{p+1}}-\overline{c_{i_{p+1}}}) \\&= \prod_{j=1}^{p+1}c_{i_j} +\Big\{ \text{lower order terms as a polynomial on } c_{i_1},\ldots, c_{i_{p+1}}\Big\} \\ & = \prod_{j=1}^{p+1}c_{i_j}
			\end{flalign*}  because again, the lower order terms are all $0$ for the exact same reason cited above.
			
			Now as for $p=2$,  we have \begin{align*}
				t\alpha_i = (e-h)(c_i-\overline{c_i}) = ec_i- c_ih+e\overline{c_i}+ h\overline{c_i} \\ =ec_i- c_ih+e\overline{c_i}
			\end{align*} because the alternating action of ${S_2}$ kills $H^0(C^2)\otimes H^*(\P E^{r+1}_n)$, and in turn, $h\overline{c_i}$.
			This give us the algebra structure on the $E_2$ page for $p\leq n-2g$ and thus completes the proof of Theorem \ref{thm3}.
			
		\end{enumerate}
		
	\end{proof}

		\section{Moduli space of smooth sections of $\grd$ on a smooth projective curve}\label{sec4}
		Le $X$ be a smooth projective curve over $\C$ of genus $g$. A line bundle $L$ on $X$ of degree $d$ is called \emph{$m$-very ample} if for every effective divisor $\xi\in X$ of degree $m+ 1$, the evaluation map
		$$ev_{\xi}: H^0(X, L) \to H^0(X, L \otimes \mathcal{O}_{\xi})$$
		is surjective, or equivalently, if $\dim H^0(X,L\otimes \mathcal{O}(-\xi)) =\dim H^0(X,L)-(m+1).$
		More generally, we have the following well-known definition. 
		\begin{defn}\label{veryampleness}
			If $\mathcal{V}$ is a \emph{linear series of type $\mathfrak{g}^r_d$} on $X$, i.e. $\V\subset H^0(X,L)$ of rank $r+1$, for some line bundle $L$ of degree $d$ on $X$, then we say $\V$ is \emph{$m$-very ample} if for every effective divisor $\xi \in X$ of degree $m+1$, we have that $$\dim\V(-\xi) = r+1-(m+1)$$ where $\V(-\xi):= H^0(X,L\otimes \mathcal{O}(-\xi))\cap \V$.
		\end{defn}
		Therefore, $0$-very ampleness is the same as global generation and $1$-very ampleness is our usual notion of very ampleness. For $m \geq 2$, the $m$-very ampleness of $\V$ is equivalent to saying that the image of $X$ under the embedding induced by $\V$ i.e. \begin{align*}
			\phi_{\V}: X\hookrightarrow  \P(\V^*)\\ x\mapsto [s_0(x):\ldots: s_r(x)]
		\end{align*} (where $s_0,\ldots,s_r$ is a basis of $\V$ as a $\C$-vector space), has no
		$(m+1)$-secant $(m-1)$-plane (note that the existence of an $(m+1)$-secant $(m-1)$-plane is special, because the expected dimension of a $(m+1)$-secant plane $m$, i.e. the span of $(m+1)$ points on $\phi_{\V}(X)\subset \P(\V^*)$ is $m$ for a general set of $m+1$ points; a wonderful reference for this is \cite[Chapter VIII]{ACGH}).
		
		In this section we are interested in the (stable) cohomology of the moduli space of smooth sections of an $m$-very ample $\grd$. We will soon see that the stability comes from the `extent' or \emph{degree} of very ampleness  (see Lemma \ref{vbness}).\begin{remark}
			Unsurprisingly, there is no necessary and sufficient condition for a linear system to be $m$-very ample that is solely determined by the parameters $g, r$ and $d$. However, there are various estimates on $m$, some of which give necessary, and some sufficient conditions for when a $\grd$ is $m$-very ample.
			\begin{itemize}\item For sufficient conditions, Farkas, in \cite{Farkas2008}, says that given a general genus $g$ smooth projective curve $X$, if we have the following inequality $$ \rho(g,r,d) -(r-m+2)+m\leq 0$$ where $\rho(g,r,d):= g- (r+1)(g-d+r)$ is the Brill-Noether number (see \cite[Section 1, Chapter IV]{ACGH}) then there exists a $\grd$ that is $m$-very ample.
				
				\item In the same paper, Farkas states a series of inequalities in Theorem 0.5, which, when simultaneously satisfied, provide sufficient  conditions for the existence of a $\grd$ that is not $m$-very ample.
				\item  In \cite[Chapter VIII]{ACGH}, for $\V \subset H^0(X,L)$  a $\grd$, they compute the virtual fundamental class of the degeneracy loci of the evaluation map of the following vector bundles on $\Sym^m X$:
				
				\begin{center}
					\begin{tikzcd}[column sep=small]
						\V\times \Sym^m X \arrow[dr]  \arrow{rr}{ev} && E_L \arrow[dl]\\
						&\Sym^m X 
					\end{tikzcd} 
				\end{center}
				where the stalks of the vector bundle $E_L$ is at a point $\xi \in \Sym^mX$ is given by $$(E_L)_\xi = H^0(X,L/L(-\xi)).$$ Note that the degeneracy loci being empty corresponds to $\V$ being $m$-very ample. In theory one can deduce inequalities involving $g,r,d$ and $m$ for which the virtual fundamental class is empty, as Farkas does in \cite{Farkas2008}, for most purposes, the formula is extremely complicated and unyielding.
			\end{itemize}
		\end{remark}  
		The vector space $\V\subset H^0(X,L)$ contains, as a Zariski open dense subvariety, the locus of smooth sections $\V^\circ$, i.e. $$\V^\circ :=  \{s\in \V: \nexists x \in X \text{ such that }v_x(s)\geq 2 \}$$ where $v_x(s)$ denotes the order of vanishing of $s$ at $x$.
		
		Geometrically, when $\V$ is $m$-very ample with $m\geq 1$, the image of $X$ under the induced embedding $\phi_\V: X\to \P(\V^{*})$ is a smooth projective curve of degree $d$ and up to $\C^*$ an element of $\V^{\circ}$ determines, and is determined by, a hyperplane in $\P(\V^*)$ that intersects $\phi_\V(X)$ smoothly i.e. at exactly $d$ distinct points.
		
		Our goal, now, is to compute the (stable) cohomology $H^*(\V^\circ;\Qb)$. 
		\begin{proof}We fix $\V$, an $m$-very ample $\grd$ for the rest of this section.
			\begin{enumerate}[Step 1.]
				\item 	\textit{Construction of a $\fS$-object $T_{\bullet}$.}
				
				First we construct a $\Delta S$ object in the category of schemes augmented on the `discriminant locus' $Z:= \V-\V^{\circ}$.
				Define $$T_0 := \{(s,x): s\in \V, v_{x}(s)\geq 2\}$$ i.e. $T_0$ is the normalisation of the discriminant locus $Z$; indeed $T_0$ is smooth,  \begin{align*}
					\pi_0:	T_0\to Z\\ (s,x) \mapsto s
				\end{align*} is a finite surjective morphism, and an isomorphism over a Zariski open dense subset of $T_0$ given by the locus of sections which are singular at exactly one point in $X$. Now for $p\geq 0$ define $$T_{p} := \overline{T_0^{\times_{Z}{(p+1)}} -\{\text{all diagonals}\}}$$ where $T_0^{\times_{Z}(p+1)}$ is the $(p+1)$-fold fibre product over $Z$, and for convenience that will be clear later, we set $$T_{-1}:= \V,$$ (deviating from the standard texts that define $T_{-1}$ to be $Z$.) Equivalently, for $p\geq 0$, we have
				\begin{flalign*}
					T_{p}=\Big\{(s,(x_0,\ldots, x_p)): \mathit{div}(s) \geq 2{\textstyle \sum} x_i\Big\}
				\end{flalign*} where $ \mathit{div}(s)$ denotes the divisor of $s \in \V$.
				
				Henceforth, unless otherwise mentioned, we use $T_{\bullet}$ to mean the semisimplicial space $T_{p\geq 0}$. Clearly  $T_{\bullet}\to Z$ is a symmetric semisimplicial object augmented over $Z$, with face maps corresponding to forgetting one of the factors of $X$:
				\begin{gather*}
					f^i:T_p\to T_{p-1} \\ s,(x_0,\ldots, x_{p})\,\,\, \mapsto \,\,\, s,(x_0,\ldots, \widehat{x_i},\ldots, x_{p})\end{gather*}
				which are all finite morphisms, and for all permutations $\sigma_0,\ldots, \sigma_p$ of $0,\ldots, p$ under the action of ${S}_p$, we define  \begin{gather*}
					\pi_p:= f^{\sigma_0}\circ\ldots\circ f^{\sigma_p}: T_p \to Z\\ s,(x_0,\ldots, x_{p})\mapsto s. 
				\end{gather*} For simplicity we abuse notation and denote, for all $p\geq 0$, the composition $$T_p\xrightarrow{\pi_p} Z\xhookrightarrow{\iota} \V$$ by $\pi_p$ as well, instead of $\iota_{\circ} \pi_p$. On the other hand, for all $p\geq 0$ we have the other projection map: \begin{gather*}
					\psi_p: T_p\to X^{p+1} \\ s, (x_0,\ldots, x_{p})\mapsto (x_0,\ldots, x_{p});
				\end{gather*} this will be particularly useful in the next step.
				
				\item 	\textit{Geometry of $T_p$.}
				
				For each $p\geq 0$, define a vector bundle $E_p\to X^{p+1}$ as follows (a similar construction is followed in \cite[Chapter IV, Section 2]{ACGH} over $\Sym^{p+1} X$). Let  $D(p+1) \subset X\times X^{p+1}$ be defined by $$D(p+1) := \{x,(x_0,\ldots, x_{p}): x=x_i \text{ for some } 0\leq i\leq p\},$$ and let $pr_j$, for $j=1,2$ denote, respectively, the projection to the first factor $X$ and the second factor $X^{p+1}$. Then $${E}_p: = {pr_2}_*(\mathcal{O}_{2{D(p+1)}}\otimes pr_1^*L).$$ is a locally free sheaf (because, $pr_1^*L$ is locally free and $2D(p+1)$ being flat over $X^{p+1}$ imply $\mathcal{O}_{2D(p+1)}\otimes pr_1^*L$ is also flat) and equivalently a vector bundle, with stalks given by $$\big({E_{p}}\big)_{\overline{x}} = H^0(X, L\otimes \mathcal{O}_{2\xi(\overline{x})})$$ where for all $\overline{x}\in X^{p+1}$ we define $\xi(\overline{x})$ to be the corresponding unordered $(p+1)$-tuple. We will often abuse notation and denote an divisor by $\xi$ when there is no scope of confusion.
				
				The natural map of sheaves given by restriction $$pr_1^*L\to \mathcal{O}_{2D(p+1)}\otimes pr_1^*L$$ induces a map on pushforwards called the evaluation map $$ev: {pr_2}_* pr_1^*L = H^0(X,L)\otimes \mathcal{O}_{X^{p+1}} \to E_p,$$ where $H^0(X,L)\otimes \mathcal{O}_{X^{p+1}}$ is the trivial bundle on $X^{p+1}$ with fibres $H^0(X,L)$. Restricting this to  $\V\subset H^0(X,L)$ gives us the map $$ev: \V\otimes \mathcal{O}_{X^{p+1}} \to E_{p+1}.$$
				At the level of stalks, one obtains the map $$\V\to H^0(X,L\otimes \mathcal{O}_{2\xi(\overline{x})}),$$ the kernel of which is precisely $\V(-2\xi(\overline{x}))$. Indeed, for any divisor $\xi$ on $X$, we have a short exact sequence of locally free sheaves on $X$:
				$$0\to L(-2\xi)\to L \to L\otimes \mathcal{O}_{2\xi} \to 0$$ that induces a long exact sequence of cohomology $$0\to H^0(X,L(-2\xi))\to H^0(X,L)\to H^0(X,L\otimes \mathcal{O}_{2\xi}) \to \ldots$$ and taking intersection with $\V$ gives us that at the level of stalks $$\mathrm{Kernel}\Big((ev_{\overline{x}}): \V\to  H^0(X,L\otimes \mathcal{O}_{2\xi(\overline{x})})\Big) = \V(-2\xi(\overline{x})).$$
				And now note that we have the following diagram:
				\begin{center}
					\begin{tikzcd}[column sep=scriptsize]
						T_{p} \arrow{rr} \arrow{drr}{\psi_p}	&&\V\times X^{p+1}  \arrow[d]  \arrow{rr}{ev} && E_p \arrow{dll}\\
						&& X^{p+1} 
					\end{tikzcd} 
				\end{center}
				where, by definition, $\psi_p: T_{p}\to X^{p+1}$ is the kernel of the evaluation map. And now observe that by the very definition of $m$-very ampleness (see Definition \ref{veryampleness}), we obtain the following:
				\begin{lemma}[Stable bound for cohomology]\label{vbness}
					For each $0\leq p\leq {m+1\over 2}$, we have \begin{align*}
						\psi_p: T_p\to X^{p+1} \\ 
					\end{align*} is a vector bundle with the fibre over a point $(x_0,\ldots,x_p)\in X^{p+1}$ given by $$\psi^{-1}(x_0,\ldots, x_p) = \V(-2\xi(x_0, \ldots, x_p)) \cong \C^{r+1-2(p+1)}.$$
				\end{lemma}
				
				\item 	\textit{Constructing a spectral sequence for the semisimplicial object $T_\bullet$.}
				
				Let $j:\V^{\circ}\hookrightarrow \V$ denote the inclusion. Then by Lemma \ref{lem2.11} we have an acylic complex of sheaves of $\Qb$-vector spaces on $\V$ given by: \begin{align*}
					j_!\Qb_{\V^{\circ}}\to \Qb_{\V}\to {\pi_0}_*\Qb_{T_0}\to \Big({\pi_1}_*\Qb_{T_1} \otimes \sgn_2\Big)^{{S}_2}\to \ldots \to \big({\pi_{p}}_*\Qb_{T_{p}}\otimes \sgn_{p+1}\big)^{{S}_{p+1}} \to \ldots
				\end{align*}
				Let $C^{\bullet}$ denote the complex $$\Qb_{\V}\to {\pi_0}_*\Qb_{T_0}\to \Big({\pi_1}_*\Qb_{T_1} \otimes \sgn_2\Big)^{{S}_2}\to \ldots. $$ Let us define $T_{-1}=\V$, and $\pi_{-1}:= id_{T_{-1}}$, the identity map on $T_{-1}$.
				Taking $\mathbf{R}\Gamma_c$ of the complex $C^{\bullet}$ one obtains a spectral sequence that reads as \begin{align}\label{sseqX}
					E_1^{p,q}= R^q\Gamma_c\Big(\V, \big({\pi_{p-1}}_*\Qb_{T_{p-1}}\big)^{{S}_p} \Big) \implies R^{p+q}\Gamma_c(\V, j_!\Qb_{\V^{\circ}})
				\end{align}
				On the right hand side, we have $$R^{p+q}\Gamma_c(\V, j_!\Qb_{\V^{\circ}}) = H_c^{p+q}(\V_{\circ};\Qb).$$ To simplify the $E_1^{p,q}$ terms we go through the following steps:
				\begin{flalign*}
					R^q\Gamma_c\Big(\V, \big({\pi_{p-1}}_*\Qb_{T_{p-1}}\big)^{{S}_p} \Big) \\ =\Big(H_c^q(T_{p-1})\otimes \sgn_{p}\Big)^{{S}_{p}} (-p) && (\pi_p \text{ finite}) \nonumber \\
					\cong \Big(H_c^q(X^{p}\times \C^{r+1-2p})\otimes \sgn_{p}\Big)^{{S}_p} (-p) && (\text{for all } p\leq {m\over 2,}\text{ by Lemma \ref{vbness}}) 
				\end{flalign*}
				\begin{equation*}
					\cong \begin{cases}
						H^2_c(X)\otimes \Sym^{p-1}H^1_c(X) (-p),& q=2(r+1) -3p+1\\ 
						H^0_c(X)\otimes H^2_c(X)\otimes \Sym^{p-2}H^1_c(X) (-p) \bigoplus \Sym^p H^1_c(X) (-p),  & q= 2(r+1)-3p\\
						H^0_c(X)\otimes \Sym^{p-1}H^1_c(X) (-p), & q= 2(r+1)-3p-1 \\ 0, & q \text{ otherwise }
					\end{cases}\label{Xweights} 
				\end{equation*} where the last step comes from the Macdonald's result on the permutation action of the symmetric group $S_p$ (twisted by the sign representation) on the cohomology $H^*(X)^{\otimes p}$ (see \cite{Macdonald1962}).
			
			Now observe the following: from Lemma \ref{lem2.11} we have an equality $\mathrm{R}\Gamma_c(\mathcal{V},C^{\bullet}) = \mathrm{R}\Gamma_c (\mathcal{V}, j_{!}\Qb_{\mathcal{V}}^{\circ})$ in the derived category of constructible sheaves over $\mathcal{V}$; on the other hand, for any $N\in \mathbb{N}$ we have $$\mathrm{R}^i\Gamma_c(\mathcal{V}, C^{\bullet}) \cong \mathrm{R}^i\Gamma_c(\mathcal{V}, C^{\bullet}/\tau_{\geq N} C^{\bullet})$$ for all $i\geq 2(r+1)-2N$, where $\tau_{\geq N} C^{\bullet}$ denotes the truncated complex up to the $(N-1)$ term and this is because $\tau_{\geq N} C^{\bullet}$ is supported on complex codimension $N$ in $\mathcal{V}$.
			
			This observation, paired with Poincaré duality $H^i_c(\V^{\circ};\Qb) \cong H^{2(r+1)-i}(\V^{\circ};\Qb)$ gives us that for all $i\leq {m-1\over 2}$: \begin{align}
				H^i(\V^{\circ};\Qb) \cong \begin{cases}
					\Sym^{p-2}H^1(X;\Qb)(-(p-1)) \oplus \Sym^p H^1(X;\Qb)(-p) & i = 2p\\
					\Sym^{p-1} H^1(X;\Qb) (-(p-1)) \oplus \Sym^p H^1(X;\Qb)(-(p+1)) & i=2p+1
				\end{cases}
			\end{align} which completes the proof of Theorem \ref{thm5}.		
			
		\end{enumerate}
		
	\end{proof}

	\begin{remark}
		Note that we could have defined the $\Delta S$-space $\tilde{T}_{\bullet}$ by defining $\tilde{T}_p:= T_0^{\times_{Z}(p+1)}$ and directly used Lemma \ref{lem2.8}. The rest of the computation, and in particular the analysis of the spectral sequence would have remained exactly the same though, because the permutation action of the symmetric groups twisted by sign kills all degeneracies.
	\end{remark}

	\bibliographystyle{alpha}
	\bibliography{SemisimplicialMathZ}
	
\end{document}